\documentclass[11pt]{amsart}

\usepackage{graphicx} 
\usepackage[margin=1in]{geometry}
\usepackage{amsmath,amssymb,amsthm}
\usepackage{bm,color,hyperref}
\usepackage{caption}
\usepackage{subcaption}
\usepackage{booktabs}

\newtheorem{theorem}{Theorem}

\newtheorem{remark}[theorem]{Remark}
\newtheorem{example}[theorem]{Example}
\newcommand\bC{{\mathbb C}}
\newcommand\bR{{\mathbb R}}
\newcommand\bZ{{\mathbb Z}}
\newcommand{\sO}{{\mathcal O}}
\newcommand{\hatz}{{\hat z}}
\newcommand{\zhat}{\hatz}
\newcommand{\dz}{\delta z}

\newcommand{\zest}{{\tilde z}}

\newcommand{\nuMax}{\nu_{\rm max}}
\newcommand{\sJ}{{\mathcal J}}

\title[Spiral endgames]{Spiral endgames \\ for computing singular endpoints of homotopy paths}                      
\begin{document}

\author{Emma R. Cobian}
\address{Rose-Hulman Institute of Technology, Terre Haute, IN, USA}
\email{cobian@rose-hulman.edu}
\urladdr{www.emmacobian.com}

\author{Jonathan D. Hauenstein}
\address{University of Notre Dame, Notre Dame, IN, USA}
\email{hauenstein@nd.edu}
\urladdr{www.hauenstein.phd}

\author{Charles W. Wampler}
\address{University of Notre Dame, Notre Dame, IN, USA}
\email{cwample1@nd.edu}
\urladdr{www.cwampler.com}

\begin{abstract}
We develop a new endgame for numerically computing singular solutions to complex analytic systems by means of a homotopy path whose endpoint is the solution to be computed. When the solution is singular, prediction-correction path tracking may fail at the endpoint. Instead of tracking the path all the way to its end, endgames sample nonsingular points on the incoming path and approximate the endpoint by interpolating the samples with a truncated Puiseux series.  By collecting sample points along a logarithmic spiral in the complex plane, the proposed spiral endgame balances between the superior numerical conditioning of a circular sample (known as the Cauchy endgame) and the update efficiency of a linear sample (known as the power series endgame), these being tied together as two extremes of the new approach.  The numerical conditioning of the methods are analyzed and compared, and computational examples illustrate the effectiveness of the new spiral 
sampling method compared with the
circular sampling method of the Cauchy endgame.


\end{abstract}

\maketitle

\section{Introduction}

Homotopy methods compute solutions to a complex analytic system of equations $f(z)=0$, 
\mbox{$f:\bC^N\rightarrow\bC^N$}, by defining continuous paths that lead from known start points to unknown solution points. This is accomplished by constructing a homotopy $h(z,t):\bC^N\times\bC\rightarrow\bC^N$ such that $h(z,0)=f(z)$ and one or more solutions to $h(z,1)=0$ are known. By varying $t$ and numerically tracking the path emanating from a known nonsingular solution at $t=1$, there is an implicitly defined smooth path $z(t)$ in an open neighborhood of $t=1$, and if that path remains nonsingular 
and bounded
for $t\in(0,1]$, the method succeeds in finding a solution to $f=0$ in $\bC^N$.

Let  $z(0)=\lim_{t\rightarrow0}z(t)\in\bC^N$ be the endpoint of
such a homotopy path. If the endpoint is nonsingular, that is, if $\det Jf(z(0))\neq0$, where $Jf(z)$ is the Jacobian matrix of $f$ evaluated at $z$, then
Newton's method refines the accuracy of the approximated endpoint with a quadratic rate of convergence. 
However, when the endpoint is singular, 
the convergence rate of Newton's method in the neighborhood slows down and may even diverge~\cite{GO83}, and the final accuracy of the estimated solution may degrade. Singular endgames~\cite{ParallelEndgame,PolyhedralEndGame,MSW91,MSW92PSEG,SWS96} remedy the situation by sampling the path in the vicinity of the endpoint where the samples are still nonsingular and 
then uses those samples to build a local model of the homotopy path. 
This local model provides a high-accuracy numerical estimate of
the endpoint~$z(0)$. 

As we shall describe, two leading singular endgames in use are the ``power series endgame''~\cite{MSW92PSEG} and the ``Cauchy endgame''~\cite{MSW91}. The power series method is more efficient when the winding number
is small, but can struggle for paths which
have a high winding number.  In contrast, the Cauchy endgame requires more computation, but it has superior numerical characteristics that help it succeed in cases of high winding number. The key difference between the two schemes is the manner in which they sample the homotopy path: the power series endgame samples at real values of $t$, whereas the Cauchy endgame samples around a circle centered on the origin. The novel contribution of this article is a new sampling scheme, the {\em spiral sample}, that provides a continuous family of methods that includes both the Cauchy scheme and the power series scheme
as two extremes.  This provides the opportunity to split the difference between the two existing schemes, thereby 
improving the numerical conditioning of the power series 
method while retaining some of~its~efficiency.

Like its predecessors, the spiral endgame, 
which is introduced in \S\ref{sec:Spiral},
is applicable to any complex analytic homotopy. Using the notation that $V(f)$ is the solution set of $f=0$, i.e.,
\[
V(f)=\{z\in\bC^N~|~f(z)=0\},
\]
three common homotopy formulations are: the Newton homotopy
\[
    h(z,t)=f(z)-tf(z_1)=0,\quad z(1)=z_1,
\]
where $z_1\in\bC^N$ is a random starting guess; the parameter homotopy
\begin{equation}\label{eq:ParameterHomotopy}
h(z,t)=f(z,tp_1+(1-t)p_0)=0, \quad z(1)\in V(f(z,p_1)),
\end{equation}
where $f(z,p):\bC^N\times\bC^M\rightarrow\bC^N$ 
is a parameterized system
with $p_1\in\bC^M$ being a set of parameters where at least one nonsingular solution is known and $f(z,p_0)=0$ being the target system to solve; and the linear homotopy
\begin{equation}\label{eq:LinearHomotopy}
    h(z,t)=(1-t)f(z)+\gamma tg(z),\quad z(1)\in V(g(z))
\end{equation}
where $g:\bC^N\rightarrow\bC^N$ is a start system whose nonsingular solutions are known or are easily computed and $\gamma\in\bC$ is randomly selected. In the particular case that $f$ is a polynomial system, the methods of numerical algebraic geometry, e.g., see~\cite{NAGbook,BertiniBook,SW05},
construct a start system $g$ and its solutions with a probability-one guarantee that every path is nonsingular for $t\in(0,1]$ and at least one path leads to every isolated solution of a system of polynomial equations. 
Moreover, the number of paths leading to each isolated solution of $f(z)=0$ is 
equal to its multiplicity with respect to $f$. For polynomials, the linear homotopy, which we will use in the 
illustrative examples of this paper, is a special case of \eqref{eq:ParameterHomotopy} where the parameters are the coefficients of the polynomials.

The proposed spiral endgame shares with the Cauchy endgame an efficiency measure called backtracking, first described in~\cite{ParallelEndgame}. This is relevant when one uses a homotopy to find all isolated solutions of a polynomial system. The main cause of the extra cost of the Cauchy endgame, as compared to the power series endgame, stems from the fact that to apply the Cauchy integral formula, one must track a circle around the origin
a positive integer number of times, stopping when the associated value of $z$ equals its starting value on the circle. This number, $\omega$, is the winding number of the path, sometimes also called the cycle number. This means that in the process of collecting samples around a circle of radius $\rho$, we have computed $\omega$ different solutions to $h(z,\rho)=0$. For each of these, the backtracking procedure tracks backwards from $t=\rho$ to $t=1$, where it arrives at a start point of the homotopy. Each such start point can be marked as completed since they all share the same endpoint. This eliminates the expense of executing the endgame $\omega$ different times to obtain the same endpoint each time. As we shall see, 
except in the extreme case 
of the power series scheme,
the spiral endgame 
provides this same capability.

The rest of this article is structured as follows.
Section~\ref{sec:Puiseux} 
provides the necessary background
on Puiseux series approximation
which is the mathematical
foundation 
of the 
Cauchy and power series endgames
used in Section~\ref{sec:Background}.
Section~\ref{sec:Spiral} describes spiral sampling 
and the spiral endgame. 
Numerical conditioning of the endgames is analyzed in Section~\ref{sec:condition}.
Several examples are considered in Section~\ref{sec:Examples}.
A short conclusion is provided in Section~\ref{sec:Conclusions}.

\section{Puiseux series approximation}\label{sec:Puiseux}
All of the endgames we discuss rely on the same underlying mathematics: the endpoint of the homotopy path is approximated by interpolating samples near $t=0$ with a truncated Puiseux series. Accordingly, we begin with a 
numerical analysis of this process.

Suppose that $h(z,t):\bC^N\times\bC\rightarrow\bC^N$ is complex analytic
and $z(t):[0,1]\rightarrow\bC^N$ is a solution path 
satisfying $h(z(t),t)\equiv0$ such that $z(t)$ is smooth for $t\in(0,1]$.
The total differential of $h(z,t)=0$ yields the Davidenko 
differential equation:
\begin{equation}\label{eq:Davidenko}
    \frac{\partial h}{\partial z} \cdot \dot{z}(t) + \frac{\partial h}{\partial t}=0.
\end{equation}
Smoothness of $z(t)$ for $t\in(0,1]$ provides that
\begin{equation}\label{eq:DerivativeZ}
\dot{z}(t) = -\left(\frac{\partial h(z(t),t)}{\partial z}\right)^{-1}\cdot \frac{\partial h(z(t),t)}{\partial t}
\end{equation}
is well defined for $t\in(0,1]$. Numerical continuation methods use predictor-corrector algorithms to track solution paths. These alternate between a prediction phase based on standard ODE solvers that advances the solution in $t$, followed by several iterations of Newton's method applied to $h(z,t)=0$ with $t$ held constant to remove prediction error, 
e.g., see~\cite{allgower2012numerical,BertiniBook}.

With this setup, the goal is to accurately compute the endpoint 
\mbox{$z(0) = \lim_{t\rightarrow0} z(t)\in\bC^N$}.
By Puiseux's theorem, e.g., see \cite[Chap.~7]{Fischer2001}, 
there exists $r>0$ such that~$z(t)$ has a convergent Puiseux series
expansion for $|t|<r$, namely 
\begin{equation}\label{eq:PuiseuxExpansion}
    z(t)= \sum_{j=0}^{\infty} c_j\cdot t^{j/\omega} = c_0 + c_1\cdot t^{1/\omega} + c_2\cdot t^{2/\omega} + \cdots
\end{equation}
where $\omega\in\bZ_{>0}$ is the winding number, also called the cycle number, of the incoming path and~\mbox{$c_j\in\bC^N$}.
The largest $r>0$ for which the series converges is called the \emph{radius of 
convergence}, and the set $\{t\in\bC~|~0<|t|<r\}$ is called the \emph{endgame operating zone}.
From~\eqref{eq:PuiseuxExpansion}, one aims to accurately compute $z(0) = c_0$.  If $f(z) = h(z,0)$, then $c_0$ is a solution of $f=0$.  

\begin{example}\label{ex:IllustrativeSample}
As an illustrative example, consider the linear homotopy
\begin{equation}\label{eq:IllustrativeHomotopy}
h(z,t) = (1-t)(z-1/2)^2 + \gamma t (z^2-1) = 
(\gamma t + 1 - t)z^2 + (t - 1) z + (1-t)/4 - \gamma t
\end{equation}
where $\gamma = e^{3i/7}$ and $i=\sqrt{-1}$.  This homotopy has two paths starting at $z=\pm1$ for $t=1$ and ending at a double root $z=1/2$ for $t=0$. 
The quadratic formula
yields that the two solution paths are 
\begin{equation}\label{eq:SolnPathzpm}
z(t) = 
\frac{-B\pm\sqrt{B^2-4AC}}{2A} = 
\frac{1}{2} \pm \frac{\sqrt{3\cdot\gamma}}{2}\cdot t^{1/2} + \cdots, ~~~~
A = \gamma t + 1-t, ~~ 
B = t-1, ~~
C = (1-t)/4 - \gamma t,
\end{equation}
showing that $\omega=2$ and the endpoint is $c_0=1/2$ for both paths.  
The radius of convergence
is the minimum of the radius of convergence of $A^{-1}$ and $\sqrt{B^2-4AC}$
centered at $t=0$.  
First, 
$$A^{-1} = 
(\gamma t + 1-t)^{-1} 
= \sum_{j=0}^\infty 
(1-\gamma)^j\cdot t^j
$$
which converges whenever $|t|<|1-\gamma|^{-1} \approx 2.3513$.
Second, 
$$\sqrt{B^2-4AC}
=\sqrt{3\cdot\gamma}\cdot t^{1/2} \cdot \sqrt{1-t+\frac{4\cdot\gamma\cdot t}{3}}
= \sqrt{3\cdot\gamma}\cdot t^{1/2}\cdot \sum_{j=0}^\infty \binom{1/2}{j}\cdot t^j\cdot  \left(\frac{4\cdot \gamma}{3} - 1\right)^j
$$
which converges whenever 
$|t|<|4\gamma/3-1|^{-1} \approx 0.5935$.
In particular, the two paths
intersect at $t=(1-4\gamma/3)^{-1}$.
In summary, this shows that $r = |4\gamma/3-1|^{-1}$.  
\end{example}

For more general examples, we know that a convergent Puiseux series exists in the neighborhood of $t=0$, but to construct it exactly, one would need to already know $z(0)$, since that is the constant term of the series. We wish to work in the opposite direction by using samples of $z(t)$ near zero to approximate the Puiseux series as a means of approximating $z(0)$. 

\subsection{Sample point formulation}\label{sec:SamplePoint}

Let $\bC^\times = \bC\setminus\{0\}$
and suppose that we are given
a solution path
$z(t)$, numbers
$t_0,\dots,t_K\in\bC^\times$,
and a purported winding number 
$\nu\in\bZ_{>0}$.  
In order to use~\eqref{eq:PuiseuxExpansion}
with $\omega = \nu$,
the value of $s_k = t_k^{1/\nu}$
needs to be computed,
which {\em a priori} has $\nu$
different values.  
For any $t^*\in(0,1]$, 
we can initialize $s^* = (t^*)^{1/\nu} > 0$ as the $\nu^{\rm th}$ principal
root.
Each point $z_k = z(t_k)$ at selected values of $t=t_k$ is found by integrating the Davidenko equation \eqref{eq:Davidenko}
which also unambiguously assigns a unique value to $s_k=t_k^{1/\nu}$ for $k=0,1,\dots,K$.

When performing computations in floating-point arithmetic, one numerically integrates 
the Davidenko equation~\eqref{eq:Davidenko}
and then corrects the 
result by applying Newton's method to solve $h(z,t_k)=0$. 
This results in approximations
$\hat z_k$ of $z_k$.  Hence, 
we assume that we are given sample
points 
$(\hat z_k,t_k,s_k)\in\bC^N\times\bC^\times\times\bC^\times$ 
for $k=0,1,\dots,K$
where $s_0=t_0^{1/\nu},\dots,s_K=t_K^{1/\nu}$ are distinct.
In particular, the sampling error is 
$\delta z_k = z_k-\hatz_k$
for $k=0,1,\dots,K$.

\subsection{Endpoint approximation}\label{sec:Extrapolation}

From the sample points,
interpolation
yields an approximation
of the endpoint at $t=0$.  
Assuming that the purported winding
number $\nu$ is equal to the actual
winding number $\omega$,
the accuracy of the result is affected by two types of errors: 
sampling error, described as
$\delta z_0,\dots,\delta z_K$,
and truncation error.
Truncation 
error is the result of only
matching the first $K+1$
terms. Denoting the the contribution of the higher-order terms as $R_K(t)$,  \eqref{eq:PuiseuxExpansion}
becomes
\begin{equation}\label{eq:PuiseuxExpansionK}
    z(t)= \sum_{j=0}^{K} c_j\cdot t^{j/\omega} + R_K(t) = c_0 + c_1\cdot t^{1/\omega} + \cdots + c_K\cdot t^{K/\omega} + R_K(t) 
\end{equation}
where $R_K(t) = O(t^{(K+1)/\omega})$.
Although~\eqref{eq:PuiseuxExpansionK}
is a vector equation,
we can consider each coordinate
independently.  
Hence, we will write all equations as if $z$ were a single variable 
and apply the same procedure to each of the coordinates. 

To approximate $z(t)$ at any point $t$ inside the radius of convergence, we introduce a change coordinates $t(s)=s^\omega$:
\begin{equation}\label{eq:zts}
    z(t(s))=\sum_{j=0}^{K} c_j\cdot s^{j} + R_K(t(s))=p(s)+R_K(t(s)).
\end{equation}
Approximating $z(0)$ using $K+1$ samples at $s=s_0,\ldots,s_K$ becomes a case of polynomial interpolation. We will refer to this later as \emph{Lagrange interpolation}, which uses only points, to distinguish it from \emph{Hermite interpolation}, which uses points and derivatives.

Two classical approaches for polynomial interpolation are Lagrange polynomials, which are useful for studying the conditioning of the problem, and Newton polynomials, which lead to efficient evaluation by divided differences \cite{hamming2012numerical,kincaid2009numerical}. Using Lagrange polynomials, we have
\begin{gather}
    p(s)= \sum_{j=0}^K \ell_{j,K}(s) p(s_j)=L(s)P,\qquad\ell_{j,K}(s)= \prod_{\substack{j=0\\k\ne j}}^K\frac{s-s_k}{s_j-s_k}, \label{eq:Lagrange}\\
    L(s)=\begin{bmatrix}
        \ell_{0,K}(s)&\cdots&\ell_{K,K}(s)
    \end{bmatrix}, \qquad 
    P=\begin{bmatrix}
        p(s_0)&\cdots&p(s_K)
    \end{bmatrix}^T.
\end{gather}
The endpoint of the path is $p(0)$. However, instead of the values $p(s_j)$, we actually have samples $\zhat(s_j)=z_j+\delta z_j=p(s_j)+R_K(s_j)+\delta z_j$. Using these samples, we have the estimate
\begin{gather}\label{eq:zestLagrange}
    \zest = L(0)\zhat = L(0)(P+\delta z + R),\qquad 
    \zhat=\begin{bmatrix}
        \zhat_0&\cdots&\zhat_K
    \end{bmatrix}^T,\\
    \delta z=\begin{bmatrix}
        \delta z_0&\cdots&\delta z_K
    \end{bmatrix}^T,\qquad 
    R=\begin{bmatrix}
        R_K(t_0)&\cdots&R_K(t_K)
    \end{bmatrix}^T.
\end{gather}
Accordingly, the error in the estimate is
\begin{equation}
    \zest-z(0)=\zest-L(0)P = L(0)(\delta z + R).
\end{equation}
This gives a bound on the error of 
\begin{equation}\label{eq:absErrorBound}
    |\zest-z(0)|\le ||L(0)||\left(||\delta z|| + ||R||\right)
\end{equation}
and implies a relative condition number of
\begin{equation}\label{eq:relErrorBound}
   \kappa= \left.\frac{|\zest-z(0)|}{|\zest|}\middle/\frac{||\delta z|| + ||R||}{||\zhat||}\right.=\frac{||L(0)||\,||\zhat||}{|\zest|}.
\end{equation}

While the Lagrange form gives a convenient formula for $L$, Newton's interpolation formula based on divided differences is more efficient. It is based on the fact that since $p(s)$ is a polynomial of degree $K$ that vanishes at $K+1$ points $s=s_0,\ldots,s_K$, it can be written as
\begin{equation}\label{eq:dividedDiffFormula}
    p(s)=\sum_{j=0}^K p[s_0,\ldots,s_j]n_j(s),\qquad 
    n_j(s)=\begin{cases}
        1,&j=0;\\
        \prod_{i=0}^{j-1}(s-s_i),& j>0.
    \end{cases},
\end{equation}
where $p[s_0,\ldots,s_j]$ is a divided difference that can be computed as
\begin{equation}\label{eq:dividedDiffRecursion}
    p[s_i,\ldots,s_j]=\frac{p[s_{i+1},\ldots,s_j]-p[s_i,\ldots,s_{j-1}]}{s_j-s_i},\quad i<j,~\,\,\,~s_i\ne s_j,
\end{equation}
starting from 
\begin{equation}\label{eq:dividedDiffInitialize}
    p[s_j]=p(s_j).
\end{equation}
Using this recursion with $s=0$ in \eqref{eq:dividedDiffFormula} and the sample values
\[
p[s_j]=\zhat_j,\quad j=0,\ldots,K,
\]
one obtains the endpoint estimate $\zest$. This approach uses $\sO(K^2)$ operations.

The Puiseux series is only valid within its radius of convergence, $r$, which depends on the branch points and ramification points of the homotopy $h$.
For all but the simplest homotopies, computing these
is an exorbitant amount of work compared with just tracking the path $z(t)$.
The alternative is to use trial-and-error: compute the approximation at successively smaller radii and check the convergence properties of the sequence of approximations.

\subsection{Hermite approximation}\label{Sec:Hermite}
One can increase the order of approximation using the derivative at 
each sample point. Since $\dot{z}(t)$ given in~\eqref{eq:DerivativeZ} is already
being used in the path tracking, including this in the estimate incurs little extra computation while gaining accuracy. Neglecting the residual term $R_K$ in \eqref{eq:zts} and using the chain rule, one has
\begin{equation}\label{eq:dpds}
    p'(s)=\frac{dt}{ds}\frac{dz}{dt}=\omega s^{\omega-1}\dot{z}(t)=\frac{\omega t}{s}\dot{z}.
\end{equation}
The divided differences evaluation of $p(s)$ proceeds similarly as in \eqref{eq:dividedDiffFormula}, except that for any sample point $s_j$ where both $p(s_j)$ and $p'(s_j)$ are included, $s_j$ appears twice in the list of samples, and we add the rule
\begin{equation}
    p[s_j,s_j]=\lim_{\Delta\rightarrow 0}p[s_j,s_j+\Delta]=p'(s_j).
\end{equation}
In particular, if we include the derivative at each of $K+1$ points, the divided differences algorithm is initialized with sample points $s_0,s_0,s_1,s_2,\ldots,s_K,s_K$, and the resulting interpolating polynomial $p(s)$ matches the first $2K+2$ terms of the Puiseux series for $z(t(s))$. 

One may also generalize the Lagrange form \eqref{eq:Lagrange} to write $p(s)$ in the form
\begin{equation}\label{eq:HermiteInterpolation}
      p(s)=\sum_{k=0}^K\left[H_{k,K}(s)p(s_k) + \tilde{H}_{k,K}(s)s_kp'(s_k)\right]=H_K(s)P+\tilde{H}_K(s)\tilde{P},
\end{equation}
where the rightmost expression is just matrix notation for the sum to its left. (Note that, by \eqref{eq:dpds}, the $k^{\rm th}$ entry of $\tilde{P}$ is $s_kp'(s_k)=\omega t_k\dot{z_k}$.) Conditioning of the endpoint estimate is then governed by the norms $||H(0)||$ and $||\tilde{H}(0)||$. Symbolic formulas for these are available, but for our purposes, numeric values suffice. The entries of the matrices can be found by executing the divided differences algorithm multiple times, with a single entry in $P$ or $P'$ equal to 1 and all other entries equal to~0.

\subsection{Circular sampling and Cauchy's Theorem}\label{sec:circular}
Endpoint approximation takes on especially simple forms when the samples are collected uniformly around a complete circle in the $s$-plane, where $t=s^\omega$. This means we have samples $(z(t_k),\dot{z}(t_k),t_k,s_k)$ for $s_k=\rho^{1/\omega}e^{2\pi i k/(K+1)}$. In normalized form, this gives $\sigma_k=s_k/s_0=u^k$, where $u=e^{2\pi i/(K+1)}$ is a $(K+1)^{\rm st}$ root of unity.

For polynomial interpolation (without derivatives) on a uniformly sampled circle, 
the Lagrange form \eqref{eq:Lagrange} gives
\begin{equation}
    \ell_{j,K}(0)= \prod_{\substack{j=0\\{j\ne k}}}^K \frac{u^k}{u^k-u^j}= \prod_{\substack{j=0\\{j\ne k}}}^K \frac{1}{1-u^{(j-k)}}.
\end{equation}
Using $x^{K+1}-1=\prod_{k=0}^K(x-u^k)=(x-1)(x^K+\cdots+x+1)$, 
it follows that $\ell_{j,K}(0)=1/(K+1)$ for all $j=0,\ldots,K$. In other words, the endpoint estimate is 
just the mean of the samples:
\begin{equation}\label{eq:trapezoidal}
    \zest=\frac{1}{K+1}\sum_{k=0}^K \zhat_k.
\end{equation}

This result can be understood another way. With the substitution 
\begin{equation}\label{eq:Zsubs}
Z(s) = z(s^\omega) = \sum_{j=0}^{\infty} c_j\cdot s^{j},
\end{equation}
the Puiseux series \eqref{eq:PuiseuxExpansion} becomes a Taylor series for $Z(s)=z(s^\omega)$ that defines a holomorphic function inside the endgame operating zone. 
For $\Gamma = \{\rho^{1/\omega}e^{2\pi i \alpha}~|~\alpha\in [0,1]\}$
taken counterclockwise, Cauchy's Integral Theorem gives 
\begin{equation}\label{eq:CauchyIntegral}
  z(0)=Z(0)=\frac{1}{2\pi i}\oint_{\Gamma} \frac{Z(s)}{s}ds 
  = \int_0^{1} Z(\rho^{1/\omega} e^{2\pi i \alpha}) d\alpha 
  = \int_0^{1} z(\rho e^{2\pi i \omega \alpha}) d\alpha,
\end{equation}
that is, $z(0)$ is just the mean value of $z$ around the closed loop.
Due to periodicity, the trapezoid rule for integration 
applied to equally spaced sample points yields~\eqref{eq:trapezoidal}.
As noted in \cite{trefethen2014exponentially}, the formula is exponentially convergent. This is clear from \eqref{eq:PuiseuxExpansionK}, since doubling the number of sample points around the circle doubles the exponent in the truncation error. Due to its close relation to Cauchy's Integral Theorem, we call estimation of the endpoint from a circular sample the \emph{Cauchy~endgame}.

A Hermite version using the first derivative at each sample is leads to the approximation
\begin{equation}\label{eq:trapezoidDerives}
\zest=\frac{1}{K+1}\sum_{k=0}^K \zhat_k -
\frac{1}{(K+1)^2}\sum_{k=0}^K \omega t_k \hat{\dot z}_k.
\end{equation}

\section{Power series and Cauchy endgames}\label{sec:Background}

To compute endpoint estimates using the methods of \S\ref{sec:Puiseux}, we must (a)~collect samples along some path near $t=0$, (b)~determine a purported winding number $\nu$, and (c)~test convergence by updating the approximation at smaller radii. 
We will use a geometric sequence
of decreasing radii with ratio $\lambda\in(0,1)$.
The power series and Cauchy endgames differ in how these steps are accomplished. Both methods start by tracking the
path $z(t)$ from $t=1$ to $t=\rho$
for some $\rho\in(0,1]$ where the endgame is initiated.

\begin{figure}[!b]
    \centering
    \begin{subfigure}[m]{0.47\linewidth}
        \includegraphics[width=\linewidth]{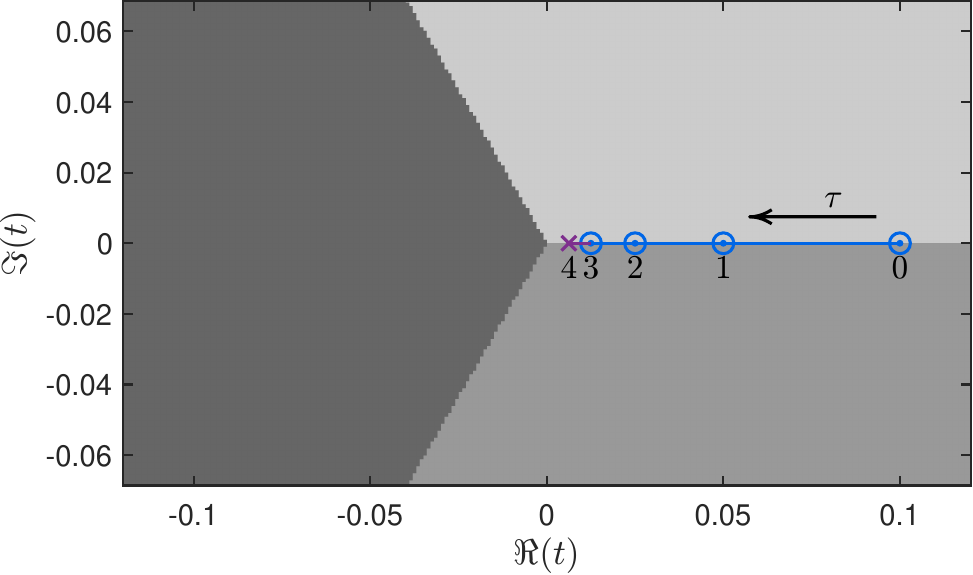}
        \caption{$t$-plane}
    \end{subfigure}
    \hspace{0.04\linewidth}
    \begin{subfigure}[m]{0.47\linewidth}
        \includegraphics[width=\linewidth]{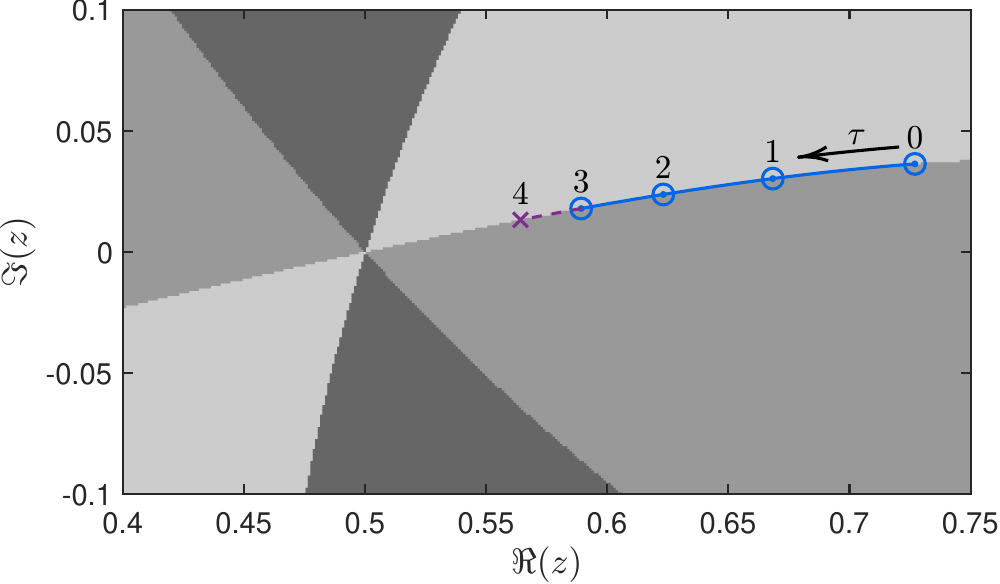}
        \caption{$z$-plane}
    \end{subfigure}
    \caption{The power series endgame applied to Example~\ref{ex:IllustrativeSample} using $\rho=0.1$, $\lambda=0.5$, $K=2$. The method tracks the homotopy path along the real line in the $t$-plane, collecting samples at a geometric series of points, $t=\rho\lambda^\tau$, $\tau\in\{0,1,2,3\}$. The resulting path and samples in the $z$-plane are also shown. Shading in both plots indicates the angle of $t$ divided into 3 wedges.  To update the estimate, the path is extended to  $\tau=4$, and the new estimate uses samples $\tau\in\{1,2,3,4\}$, which can be re-indexed as $t=\rho'\lambda^{\tau'}$, $\rho'=\rho\lambda$, $\tau'\in\{0,1,2,3\}$.}
    \label{fig:PowerEndgame}
\end{figure}

The power series endgame takes a straightforward approach: continue sampling points along the real line in the $t$-plane, i.e., along $t\in(0,\rho]$, and use the Puiseux series \eqref{eq:PuiseuxExpansion} to extrapolate to $t=0$. Specifically, its steps are as follows.
\begin{description}
    \item[Sample] Collect samples $z_0,\ldots,z_{K+1}$ at $t_k=\rho\lambda^k$, $k=0,\ldots,K+1$.
    \item[Winding Number] For trial values $\hat\nu=1,\ldots,\nuMax$, apply divided differences to the sample data at $t_0,\ldots,t_K$ to estimate $z_{K+1}$ at $t_{K+1}$ as $\breve{z}_{K+1,\hat\nu}$. Set $\nu$ to be the 
    value of $\hat\nu$ that gives the smallest approximation error $||z_{K+1}-\breve{z}_{K+1,\hat\nu}||$.
    \item[Estimate endpoint] Use the purported winding number $\nu$ and apply divided differences to the sample data at $t_0,\ldots,t_{K+1}$ to estimate $\zest$ at $t=0$.
    \item[Update] Track from $t=\rho\lambda^{K+1}$ to $t=\rho\lambda^{K+2}$. Repeat the winding number and estimation steps with $\rho\leftarrow\rho\lambda$.
\end{description}
Figure~\ref{fig:PowerEndgame} illustrates these steps being applied to Example~\ref{ex:IllustrativeSample}.

\begin{figure}
    \centering
    \begin{subfigure}[b]{0.48\linewidth}
        \includegraphics[width=\linewidth]{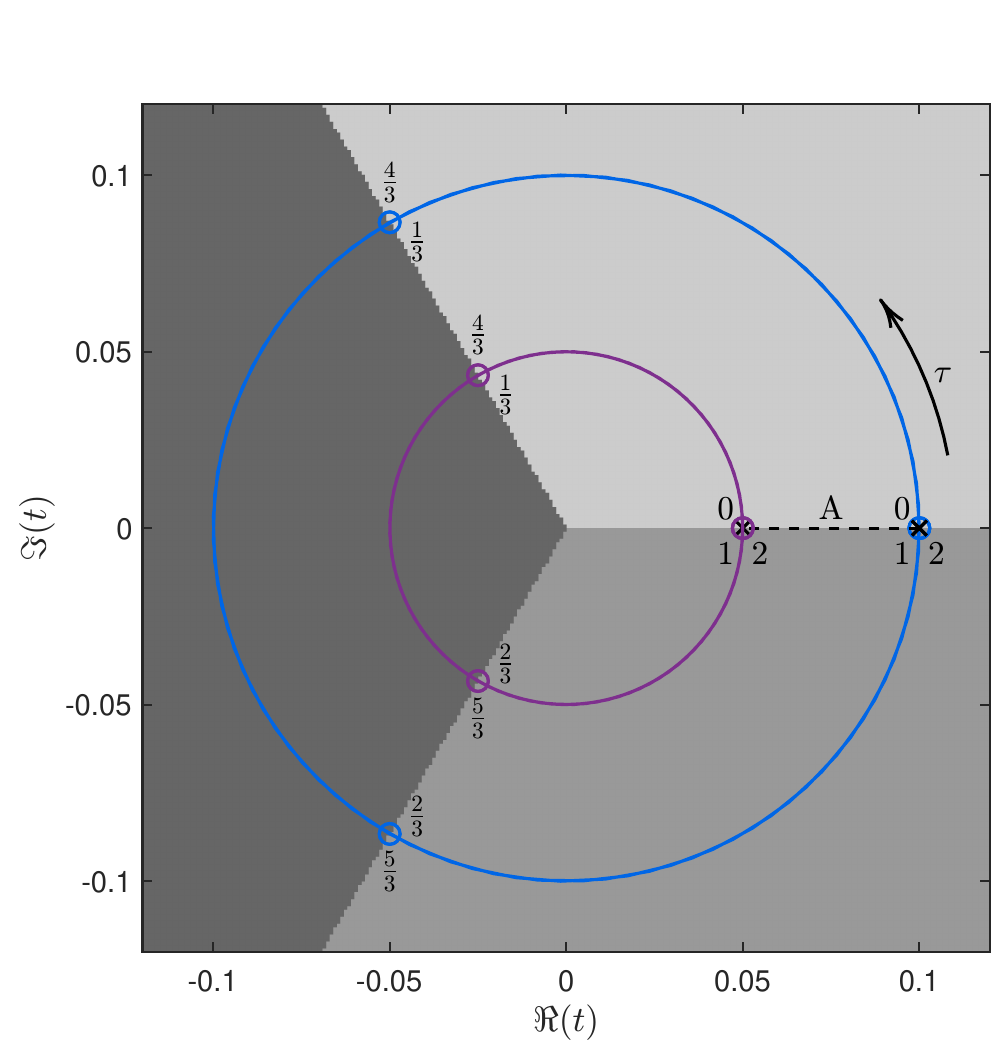}
        \label{fig:tPlane_Cauchy}
        \caption{$t$-plane}
    \end{subfigure}
    \hspace{0.03\linewidth}
    \begin{subfigure}[b]{0.47\linewidth}
        \includegraphics[width=\linewidth]{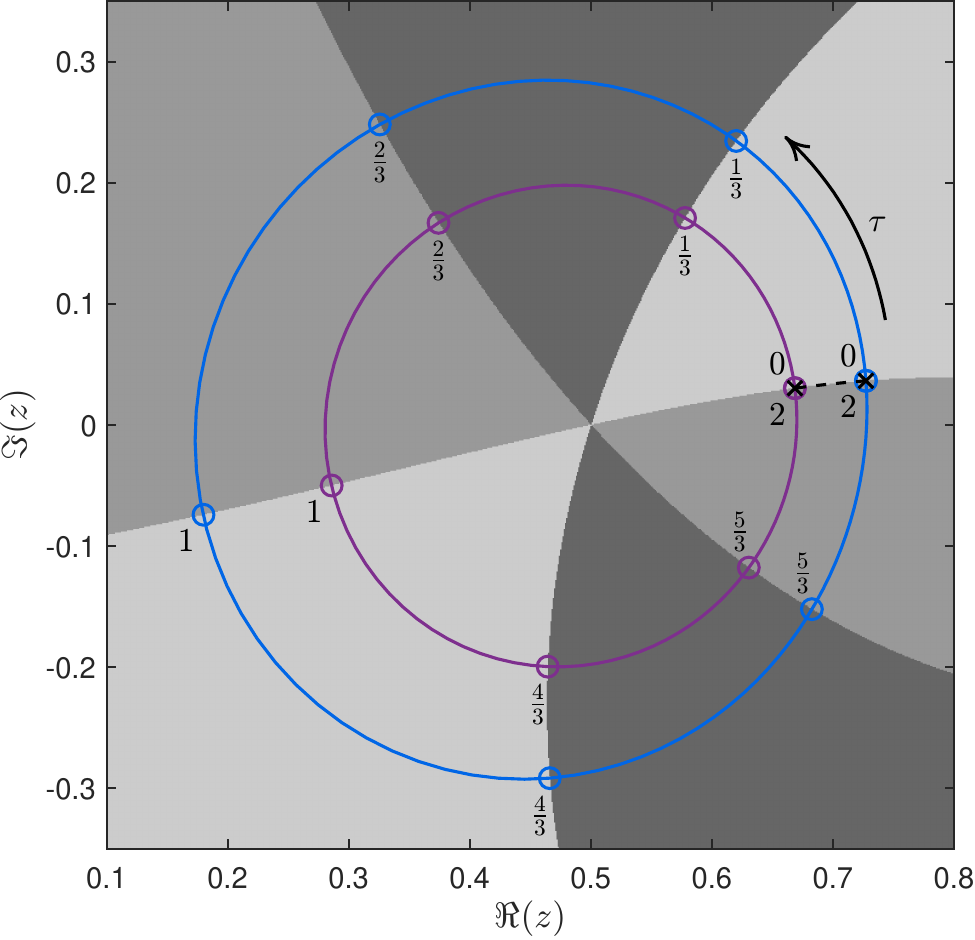}
        \label{fig:zPlane_Cachy}
        \caption{$z$-plane}
    \end{subfigure}
    \caption{The Cauchy endgame applied to Example~\ref{ex:IllustrativeSample}, with
    $\rho=0.1$, $n=3$, and $\lambda=0.5$, tracks the homotopy path around a circle in the $t$-plane, collecting samples at 3 uniformly spaced points per orbit, marked 'o'. The points are annotated with the values of $\tau$ for $t=\rho e^{2\pi i \tau}$. The corresponding path and samples in the $z$-plane are also shown. Shading in both plots indicates the angle of $t$ divided into 3 wedges. It takes two trips around the circle in $t$ to complete one revolution in~$z$. To update the estimate, the path of real $t$ is tracked along segment $A$ from $t=\rho$ to $t=\rho\lambda=\rho'$, at which point a new circular sample is initiated, again orbiting twice around the origin in~$t$ via $t=\rho'e^{2\pi i\tau}$.}
    \label{fig:CauchyEndgame}
\end{figure}

In contrast, the Cauchy endgame moves $t$ into the complex plane. By encircling the origin, it effectively shifts the use of the Puiseux series from extrapolation to interpolation as follows.
\begin{description}
    \item[Sample] Track $z(t)$ around a circular path $t=\rho e^{2\pi i\tau}$ for $\tau\in[0,\nu]$ with $\nu=1,2,\ldots$ where $i=\sqrt{-1}$. Samples are taken at $n$ equally spaced points around each orbit, resulting in $n\nu+1$ samples in total.
    \item[Winding number] Note from \eqref{eq:PuiseuxExpansion} that 
    $z(\rho e^{2\pi i\omega})=z(\rho)$ if $\rho < r$. Hence, we stop tracking the circular path at the first value of $\nu$ such that the path closes: $z(\rho e^{2\pi i\nu})=z(\rho)$. This value of $\nu$ is the purported winding number. 
    \item[Estimate endpoint]
    Since the first and last samples are equal at closure, we have $n\nu$ independent points and we apply \eqref{eq:trapezoidal} or \eqref{eq:trapezoidDerives} with $K=n\nu-1$.
    \item[Update] Track from $t=\rho$ to $t=\rho\lambda$. Repeat the sampling, winding number, and estimation steps at the new radius $\rho\leftarrow\rho\lambda$.
\end{description}
Notice that, unlike the power series approach, updating the Cauchy endgame requires repeating the sampling step anew.
Figure~\ref{fig:CauchyEndgame} illustrates these steps as they apply to Example~\ref{ex:IllustrativeSample}.

When the Cauchy endgame is applied to a polynomial system, the system $h(z,\rho)=0$ has a finite number of roots. Hence, it is clear that sampling must terminate with a closed path. For arbitrary complex analytic systems, the circular path must close when $\rho$ is small enough, but it may not do so for $\rho$ outside the radius of convergence. In that case, one should ensure termination of the 
winding number determination
step by setting a maximum value of $\nu$, say $\nuMax$.  If this limit is reached before closure, the procedure moves on to the next trial radius without producing an estimate. If the true winding number is greater than the limit, the method will keep updating at smaller and smaller radii until tracking fails or a limit on the minimum radius is reached.

These methods have advantages and disadvantages. The strongest advantage of the power series approach is that its update step only requires one additional sample. In contrast, the Cauchy update collects $n\nu+1$ samples. Since this is the most computationally expensive step of the procedures, the power series method tends to  run faster than the Cauchy method. However, the Cauchy method has several advantages of its own. First, its method of determining the winding number by detecting loop closure is more definitive than the trial-and-error test of the power series. Second, the numerical conditioning of its estimate is governed by $||L(0)||_1=1$, so it avoids the amplification of error that can arise in the power series approach. (See \S\ref{sec:condition} below for an analysis of conditioning.)
Third, in implementations where $K$ is prespecified, the power series method collects the same number of samples no matter the purported winding number, while the Cauchy method automatically increases the number of samples to $K+1=n\nu$. This means that the order of approximation of the Cauchy method is constant for all winding numbers while the power series order decreases proportional to 
the actual winding number~$\omega$. In the power series, one could respond to $\nu>1$ by collecting proportionally more samples, but this means tracking to tiny values of $t$ where $h(z,t)$ is more singular, causing tracking using \eqref{eq:DerivativeZ} to become numerically ill-conditioned.
Such ill-conditioning can be overcome using multiple precision arithmetic~\cite{AMP}, but this comes at a high cost compared to double precision.

\subsection{Winning the endgame}\label{sec:Winning}

An endgame needs to terminate when  
the endpoint estimation is deemed accurate enough. We use two successive endpoint approximations for this purpose. After an update,
the sample farthest from the origin changes from $\rho>0$ to $\rho\lambda$, $\lambda\in(0,1)$. If $\rho<r$, that is, if all the samples are within the convergence radius and $\nuMax\ge\omega$, then the winding number determined by both approximations should be equal to the actual winding number $\omega$. (Results may vary outside the endgame operating zone.) As the updates approach the origin of $t$, the endpoint approximations should converge and they should satisfy the target system $f(z)=0$.
Accordingly, we use the following stopping criteria, all of which must be satisfied.
\begin{enumerate}
    \item Successive purported winding numbers are equal.
    \item Successive endpoint approximations agree within a desired tolerance.
    \item The backward error of the final approximation is smaller than a desired tolerance.
\end{enumerate}
In Item~3, one could use the function residual, $||f(\zest)||$, instead of backward error, but this is directly dependent on the scaling of $f$. 
For a polynomial system, a criterion less sensitive to scaling is the the backward error criterion from~\cite[\S4]{EigenvalueSolver}.
In particular, if
$f_j(z) = \sum_{a} b_{j,a} z^a$
for $j=1,\dots,N$,
then the backward error is
\begin{equation}\label{eq:BWE}
{\rm BWE}(z) = \frac{1}{N}
\sum_{j=1}^N \frac{|f_j(z)|}{
1 + \sum_a |b_{j,a} z^a|}.
\end{equation}
This is the criterion we use in the examples of this paper.

\subsection{Pruning the endgame via backtracking}\label{sec:pruning}
When $f(z)$ is polynomial, one may be interested in computing all isolated solutions. Methods from numerical algebraic geometry define a homotopy and a set of start points such that the set of paths emanating from the start points includes all isolated solutions of $h(z,t)$ for each $t\in(0,1]$. 
Each isolated solution of $f(z)=0$ of multiplicity $m\geq1$ is the endpoint of $m$
paths.  The $m$ paths can be partitioned
based on the local irreducible
decomposition determined by $h(z,t)$.
Each corresponding locally 
irreducible component consisting
of $\omega$ paths will 
have a Puiseux series~\eqref{eq:PuiseuxExpansion}
approaching the isolated singularity.
If we have computed the endpoint for one of these, then there is no need to re-compute it for the others. The catch is that we do not know at the outset which start points will go to the same endpoint. 

In the endgame operating zone
it was noted in~\cite{ParallelEndgame}
that as the Cauchy endgame tracks its circular path $\omega$ times, it produces all $\omega$ solutions at $t=\rho$ on the local 
irreducible component. 
Accordingly, one may track each of these backwards from $t=\rho$ to $t=1$ to land at a start point. Removing each such start point from the queue of paths to track avoids repeating the endgame $\omega-1$ more times.

To our knowledge, this idea has never been applied to the power series endgame, but it could. Once the coefficients of the Puiseux series are known, one can use $Z(s)$ 
as in \eqref{eq:Zsubs}
to estimate each of the remaining $(\omega-1)$ values of $z(\rho)$ and backtrack to their start points.

Like the Cauchy endgame, the spiral endgame 
proposed in \S\ref{sec:Spiral}
tracks a path that generates $\omega$ points on the real line of $t$, and these can be tracked back to eliminate duplications of the endgame.

\subsection{Cauchy and power series numerical examples}\label{sec:ExamplesIntro}

Figures \ref{fig:PowerEndgame} and \ref{fig:CauchyEndgame} illustrate the power series and Cauchy endgames as applied to Example~\ref{ex:IllustrativeSample} using the homotopy $h(z,t)$ in~\eqref{eq:IllustrativeHomotopy} for the path starting at $(z,t)=(1,1)$ and sampled using $\rho = 0.1, K=2, \lambda=0.5$. The corresponding sample points are listed
below to four decimal places for the power series endgame (left) 
and the Cauchy endgame (right).

$$\begin{array}{ccc}
\text{Power series}&&\text{Cauchy}\\
\begin{tabular}{c|c|c}
   $k$  & $z(\rho \lambda^k)$ & 
   $\dot{z}(\rho \lambda^k)$
   \\ \hline
   $0$ & $0.7270 + 0.0363i$ & $0.9515 + 0.0724i$\\
   $1$ & $0.6683+0.0303i$ & $1.4899+0.1914i$ \\
   $2$ & $0.6231+0.0237i$ & $2.2606+0.3630i$ \\
   $3$ & $0.5892+0.0179i$ & $3.3587+0.6061i$    
\end{tabular} 
     & ~~ & 
\begin{tabular}{c|c|c}
   $k$  & $z(\rho e^{2 \pi i k/n})$ 
   & $\dot{z}(\rho e^{2\pi i k/n})$
   \\ \hline
   $0$ & $0.7270 + 0.0363i$ & $\phantom{-}0.9515 + 0.0724i$\\
   $1$ & $0.6199 + 0.2345i$ & $\phantom{-}0.5276 - 1.1942i$\\
   $2$ & $0.3255 + 0.2481i$ & $-0.8428 - 1.4163i$ \\
   $3$ &  $0.1796 - 0.0743i$ & $-1.9082 - 0.4159i$ \\
   $4$ &   $0.4659 - 0.2919i$ & $-1.4686 + 0.6952i$\\
   $5$ & $0.6824 - 0.1523i$ & $\phantom{-}0.0127 + 1.0117i$\\
   $6$ & $0.7270 + 0.0363i$ & $\phantom{-}0.9515 + 0.0724i$\\
\end{tabular}
\end{array}
$$ 
Since $0.1<r\approx 0.5935$, all sample points are inside the radius of convergence of the Puiseux~series.

Using $\nuMax = 4$, the power series determines a purported winding number as the value of $\hat\nu$ that minimizes $|z_3-
\breve{z}_{3,\hat\nu}|$
for $\hat{\nu} = 1,\dots,4$, namely $\nu =2$ based on the following data.
$$
\begin{tabular}{c|c|c|c|c}
$\hat{\nu}$ & $1$ & $\mathbf{2}$ & $3$ & $4$  \\ \hline \rule{0pt}{0.15in}
$|z_3-\breve{z}_{3,\hat\nu}|$ &  
$7.543\cdot10^{-3}$& $\mathbf{1.124\cdot10^{-4}}$ & $1.199\cdot 10^{-3}$ & $1.378\cdot10^{-3}$ \\ 
\end{tabular} 
$$
The following uses
derivatives via \eqref{eq:HermiteInterpolation} 
which also yields a purported
winding number of $\nu = 2$.
$$
\begin{tabular}{c|c|c|c|c}
$\hat{\nu}$ & $1$ & $\mathbf{2}$ & $3$ & $4$  \\ \hline \rule{0pt}{0.15in}
$|z_3-
\breve{z}_{H,3,\hat\nu}|$ & $9.785\cdot10^{-4}$ & 
$\mathbf{2.897\cdot 10^{-8}}$ & 
$2.585\cdot10^{-6}$ & 
$1.446\cdot 10^{-6}$  \\
\end{tabular} 
$$
Once the purported winding number $\nu$ is determined,
we can use all four points 
in the extrapolation.
This yields endpoint approximations
of $0.50012-1.35\cdot10^{-4}i$
and $0.4999999866-6.06\cdot 10^{-9}i$
without and with derivatives, respectively.

The Cauchy endgame determines the purported winding number by observing when the loop first closes, which 
first happens with $z\left(\rho e^{2\pi i (0/3)}\right)=z\left(\rho e^{2\pi i (6/3)}\right)$, giving 
the purported winding number of
$\nu=2$. The endpoint 
estimate from~\eqref{eq:trapezoidal},
which does not use derivatives,
is approximately $0.500059+6.84\cdot 10^{-5}i$,
while the estimate from~\eqref{eq:trapezoidDerives},
which does use derivatives,
is approximately
$0.5000000069-4.91\cdot 10^{-10}i$.

Note that the Cauchy endgame estimate is more accurate than the four-point power series estimate even though 
the Cauchy samples are all taken at $|t|=\rho=0.1$ while the power series had to track closer to the singularity to sample the final point at $t=\rho\lambda^3=0.0125$. This is primarily due 
to the Cauchy method having collected more points, six versus four, so it has a higher order of approximation accuracy.

\section{Spiral sampling and the spiral endgame}\label{sec:Spiral}

As described in \S\ref{sec:Background}, both  the power series and Cauchy endgames use samples of the form $z_k = z(\rho a^k)$
with differing values of $a$.
The power series endgame uses $a = \lambda\in (0,1)\subset\bR$ 
while the Cauchy endgame uses
$a = e^{2 \pi i/n}$ 
so that $a\in\bC$ with $|a|=1$.
We can think of both of these in the same family
of spiral samples by allowing $a\in\bC$ such that $0 < |a| \leq 1$.
This yields a two-real-parameter set of choices: $a=\lambda e^{2 \pi i\theta} = e^{\log \lambda + 2 \pi i \theta}$ for $\lambda\in (0,1]$ and $\theta\in [0,1)$.
In particular,~$\lambda$ describes the descent towards~$0$ while $\theta$ describes the encircling around the origin.
The Cauchy endgame has $\lambda = 1$
and $\theta = 1/n$ 
while the power series endgame has $\lambda \in (0,1)$ 
and $\theta = 0$.  
By viewing sampling in this way,
the goal is to incorporate many of the advantages of
both the Cauchy and power series endgames when $\lambda,\theta\in(0,1)$
in a computationally efficient manner as summarized below. We will see that
one downside of spiraling is that path interchanging 
can occur when outside of the endgame operating zone, but this
can be detected and mitigated.

\subsection{Spiraling and updating}\label{sec:SpiralUpdate}

Given $\rho > 0$ and $\lambda\in (0,1)$, 
consider the spiral path 
starting from $t=\rho\in(0,1]$ and limiting to $t=0$ via 
\begin{equation}\label{eq:spiralPath}
    t(\tau) = \rho e^{\tau (\log \lambda + 2\pi i)} = \rho \lambda^\tau e^{2\pi i \tau}, \qquad \tau\in[0,\infty).
\end{equation}
Hence, $t(\tau)\in (0,\rho]$ if and only
if $\tau \in \bZ_{\geq 0}$, which gives the real sequence $t(k) = \rho \lambda^k$ 
for $k\in\bZ_{\geq0}$.  Therefore,
each revolution along the spiral path 
descends by a factor of $\lambda$.

If $0 < \rho < r$ where $r$ is the radius of convergence
for $z(t)$ as in \eqref{eq:PuiseuxExpansion}, 
then the entire spiral path~$t(\tau)$ for $\tau\in [0,\infty)$
is contained in the endgame operating zone.  
Hence, the winding number $\omega$ is the minimum 
$\beta\in\bZ_{>0}$ such that one obtains a closed $z$-loop 
via the path obtained by first tracking along
the spiral using $z(t(\tau))$ from $\tau = 0$ to $\tau=\beta$
and then along the real line using 
$z(t)$ from $t = \rho$ to $t = \rho\lambda^\beta$.
Just as with the Cauchy endgame,
samples are taken at $n$ equally spaced
points in $\tau$ around each orbit.

\begin{figure}[!b]
    \centering
    \begin{subfigure}[b]{0.47\linewidth}
        \includegraphics[width=\linewidth]{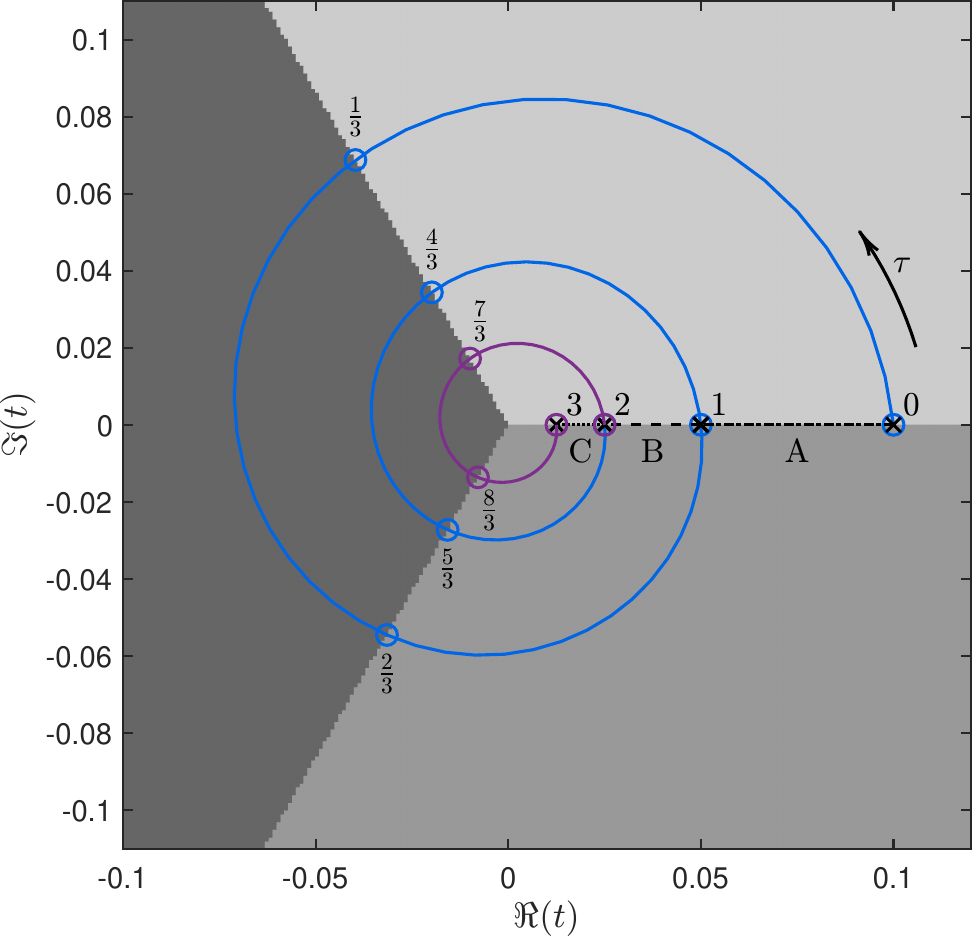}       
        \caption{$t$-plane}\label{fig:tPlane_spiral}
    \end{subfigure}
    \hspace{0.03\linewidth}
    \begin{subfigure}[b]{0.47\linewidth}
        \includegraphics[width=\linewidth]{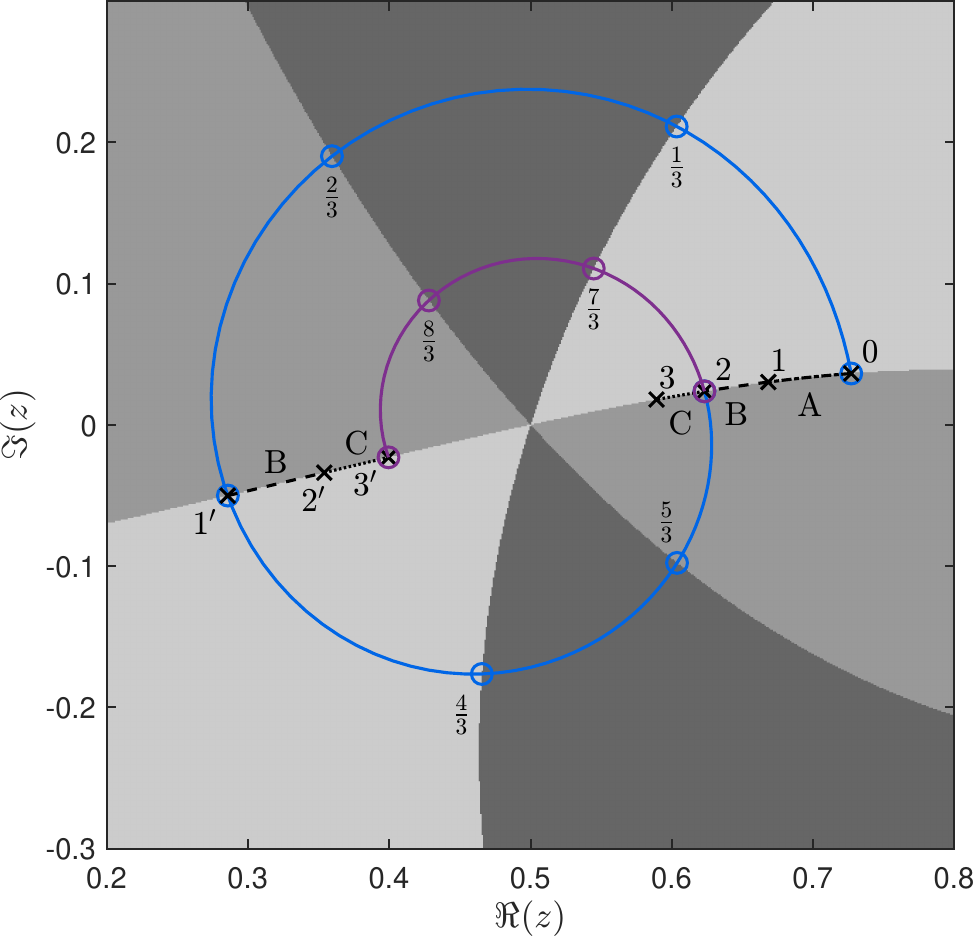}       
        \caption{$z$-plane}\label{fig:zPlane_spiral}
    \end{subfigure}
    \caption{The spiral endgame applied to Example~\ref{ex:IllustrativeSample}, with
    $\rho=0.1$, $n=3$, and \hbox{$\lambda=0.5$}, tracks the homotopy path around a logarithmic spiral in the $t$-plane, collecting samples at 3 uniformly spaced points per orbit. The points are annotated with the values of $\tau$ for $t=\rho\lambda^\tau e^{2\pi i \tau}$. Paths for real $t$ are also tracked, labeled A, B, and C. The resulting paths and samples in the $z$-plane are also shown. Shading in both plots indicates the angle of $t$ divided into 3 wedges. In the $z$-plane, spiral 0-1'-2 and arc 0-1-2 form a closed path, which is used for the first estimate of the path endpoint. To update the estimate, the spiral is tracked from $\tau=2$ to $\tau=3$, whereupon the spiral arc 1'-2-3' forms a closed path with arc 1'-2'-3'.}
    \label{fig:spiralEndgame}
\end{figure}

\begin{example}\label{ex:SpiralIllustrative}
Figure~\ref{fig:tPlane_spiral}
shows a spiral path in the $t$-plane applied to Example~\ref{ex:IllustrativeSample} with $\rho=0.1$, $\lambda = 0.5$, and $n=3$, along with a real path consisting of segments A, B, and C. The corresponding paths in the $z$-plane are shown in Figure~\ref{fig:zPlane_spiral}. After one orbit in the $t$-plane, the spiral path arrives at point~1' in the $z$-plane while tracking real segment~A arrives at point~1. Even though the $t$ values are the same, the corresponding $z$ values do not match. Hence, the procedure continues to track the spiral to $\tau=2$, where the $z$ value is found to match with the end of segment~B. This gives a purported winding number of $\nu=2$, which is equal to the actual winding number $\omega$ since point~0 started inside the endgame operating zone. 
As further illustrated in Figure~\ref{fig:spiralEndgame}, to update the approximation, one advances $\tau$ by one and finds a new closed loop. This contrasts with the Cauchy endgame, which requires $\omega$ orbits of $t$ for every update. 
The corresponding $z$ and $\dot{z}$ values for the spiral endgame are listed below to four decimal places with $\tau = k/n$. Samples along paths in real $t$ emanating from $z$ at $k=0$ and $k=3$ to determine closure are included.
\vspace{8pt}
\begin{center}
    \begin{tabular}{c|c|c|c|c}
    $k$  & $z(\rho e^{(\log \lambda + 2\pi i)k/n})$ & $\dot{z}(\rho e^{(\log \lambda + 2\pi i)k/n})$ & \textup{Real Path} $1$ & \textup{Real Path} $2$\\ 
    \hline
    $0$ & $\mathbf{0.7270 + 0.0363i}$ & $\phantom{-}0.9515 + 0.0724i$ & $\mathbf{0.7270 + 0.0363i}$ & - \\
    $1$ & $0.6035 + 0.2113i$ & $\phantom{-}0.6355 - 1.3184i$ & - & -\\
    $2$ & $0.3594 + 0.1902i$ & $-0.9634 - 1.7739i$ & - & -\\
    $3$ & $\mathbf{0.2856 - 0.0502i}$ & $-2.4241 - 0.5719i$ & $0.6683 + 0.0303i$ & $\mathbf{0.2856 - 0.0502i}$ \\
    $4$ & $0.4655 - 0.1764i$ & $-1.9728 + 1.3448i$ & - & -\\
    $5$ & $0.6038 - 0.0978i$ & $\phantom{-}0.3076 + 2.0643i$ & - & -\\
    $6$ & $\mathbf{0.6231 + 0.0237i}$ & $\phantom{-}2.2606 + 0.3630i$ & $\mathbf{0.6231 + 0.0237i}$ & $0.3540 - 0.0339i$ \\
    $7$ & $0.5448 + 0.1107i$ & $\phantom{-}1.6449 - 2.4627i$ & - & -\\
    $8$ & $0.4280 + 0.0880i$ & $-1.5033 - 3.4553i$ & - & -\\
    $9$ & $\mathbf{0.3994 - 0.0230i}$ & $-4.2746 - 1.0131i$ & - & $\mathbf{0.3994 - 0.0230i}$ \\
\end{tabular}
\end{center}
\vspace{8pt}
The endpoint estimate
at the first spiral closure 
($k=6$) is
approximately
\hbox{$0.50000117 + 6.5660\cdot 10^{-7}i$} 
without including derivatives and $0.500000000000518 - 1.4105\cdot 10^{-12}i$ with derivatives. Advancing~$\tau$ by one (ending at $k=9$) 
resulted in another spiral 
closure with the endpoint
estimate being approximately
\mbox{$0.49999985 - 5.7388\cdot 10^{-8}i$} and $0.500000000000001 - 2.4616\cdot 10^{-14}i$ without and with derivatives, respectively. See \S\ref{sec:Extrapolation} for details on approximating the endpoint.

\end{example}

One may see that once the spiral endgame is inside the endgame operating zone, each update tracks $n$ arcs to complete another orbit of the spiral plus at most $\omega$ arcs corresponding to moving along a real segment in~$t$. Hence, the update requires tracking a total of at most $n+\omega$ arcs. In contrast, the Cauchy method tracks one real arc to the next circle, and then tracks $n\omega$ arcs to complete $\omega$ orbits, for a total of $1+n\omega$ arcs.  

\subsection{Winding number determination and
path interchanging}\label{sec:interchange}

The procedure described in \S\ref{sec:SpiralUpdate} determines a purported
winding number $\nu$ by detecting closure around a loop consisting of $\nu$ orbits of the spiral path connected to one segment tracked for real $t$. We term this a ``purported'' winding number since, if the closed path starts outside the radius of convergence, we may have $\nu\ne\omega$. However, as the updates continue to spiral inward, eventually a closed path will be generated inside the radius of convergence, at which point we will have $\nu=\omega$.

If $\rho > r$, then the presence of 
a ramification point may result in a closed path that does not include any point on the arc of real $t$ emanating from the original start point. When such a \emph{path interchange} occurs, continuing the endgame from that point will give a root but not the one that is the endpoint of the path of real $t$ for the original start point. 
The following illustrates a univariate homotopy that has a path
interchange.

\begin{figure}[!t]
    \centering
    \includegraphics[width=0.9\linewidth]{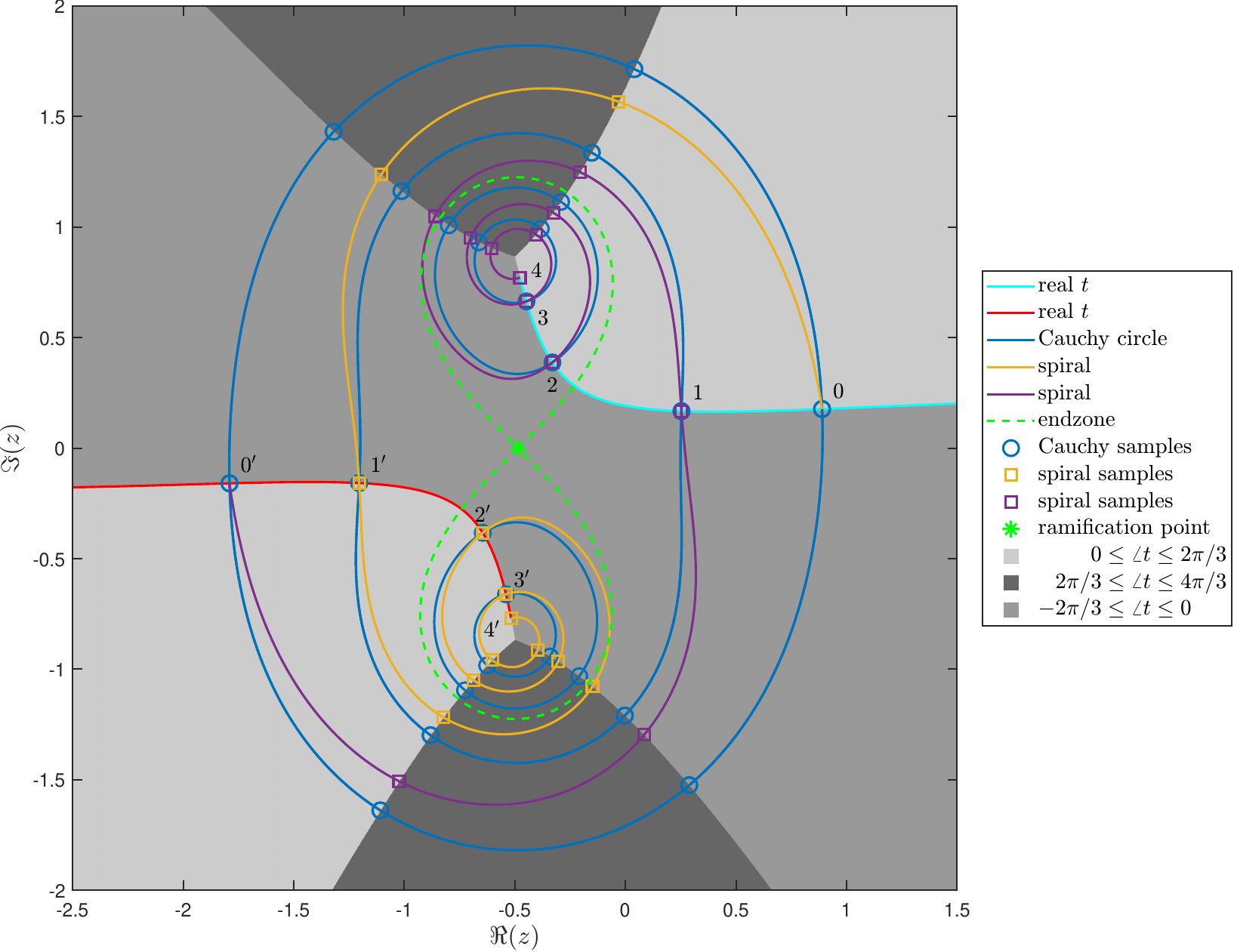}
    \caption{
    Plot of $\Re(z)$ and $\Im(z)$
    for the path interchange example 
    with homotopy~$h(z,t)$ in~\eqref{eq:pathSwap}. 
    Background shading corresponds to the angle of $t$, and roots are located where the three shades meet.} 
    \label{fig:pathswap}
\end{figure}

\begin{example}\label{ex:univariate}
    Consider the linear homotopy 
    \begin{equation}\label{eq:pathSwap}
        h(z,t)=(1-t)f(z)+\gamma t g(t)=(1-t)(z^2+z+1)+\gamma t (z^2-25),
    \end{equation}
    where $\gamma=e^{\pi i/16}$.  Figure~\ref{fig:pathswap} shows several salient features of this homotopy when the Cauchy and spiral endgames are applied. The homotopy has two ramification points, whose $z$ coordinates are the solutions to $fg'-f'g=0$. The one 
    closest to the origin in $t$ is $(z^*,t^*)\approx(-0.4853,0.0289-0.0056i)$ 
    with branch point $t^*=0.0289-0.0056i$.
    The dashed curve defined by $|t|=|t^*|\approx0.0294$ is the outer boundary of the endgame operating zone. 

    There are two homotopy paths for real $t\in[0,1]$ that start at $(z,t)=(\pm 5,1)$. These pass through the points labeled $0,\ldots,4$ and $0',\ldots,4'$
    for the positive and negative start point, respectively. 
    We begin the 
    Cauchy and spiral endgame at $t=\rho=0.1$, which corresponds to points $0$ and~$0'$, use a contraction ratio $\lambda=0.5$,
     and collect $n=3$ points per loop. The shades of gray in the image correspond to three ranges of the angle of $t$ such that the samples are collected where the endgame paths cross boundaries where two shades meet. The $z$ coordinate of the two solutions  at $t=0$ are located where all three shades meet. 
    
    First, consider the Cauchy endgame, which follows circular paths $t(\tau)=\rho\lambda^ke^{2\pi i\tau}$, for a sequence of circles $k=0,1,\ldots$. Starting from point $0$ with $k=0$ and $\tau=0$, the first trip around the circle arrives at point $0'$ when $\tau=1$. A second trip around returns to point $0$ when $\tau=2$. This gives a purported winding number of 2. The second circle travels 1,1',1, which also gives a purported winding number of 2. The Cauchy integrals for these circles, and hence the corresponding endpoint estimates, approximate the mean of the two roots inside the circles. These may pass the convergence test, but they will fail the backward error test. Circles starting at point 2 and beyond are inside the endgame operating zone. These all give 
    a purported winding number of 1 and yield increasingly accurate approximations that eventually pass both the convergence and backward error tests.

    Next, consider the spiral endgame that tracks using $t(\tau)$ from \eqref{eq:spiralPath}. We observe that the spiral starting at point 0 passes through points $0,1',2',\ldots$ while the spiral starting at $0'$ passes through $0',1,2,\ldots$. That is, the spirals have interchanged endpoints.  
    If one does not care about matching endpoints to incoming paths, and if all the endgames use the same spiral $t(\tau)$, then it is acceptable to let the paths interchange in this fashion 
    as one still obtains all roots. However, if the endgame is required to return the endpoint corresponding to the real path of $t$ from which it is initiated, it must detect the interchange. In the case at hand, the spiral $0,1',2'$ detects that the closure at the spiral sequence $1',2'$ does not involve points $1,2$ collected on real $t$. This signals that it has interchanged paths. The remedy is to restart on the spiral path emanating from point $2$. This spiral approaches the correct root, giving a constant purported winding number of 1, and progressively more accurate endpoint estimates. Starting the spiral endgame from point $0'$ gives a similar behavior.
\end{example}

Just as this example demonstrates, 
the way to identify that a path interchange
occurred is that 
there are $\alpha,\xi\in\bZ_{>0}$ such
that the following hold:
\begin{itemize}
    \item starting from the 
    same point $z(t(\alpha))$, 
    the point $z(t(\alpha+\zeta))$ obtained by tracking along the spiral
    is the same point 
    as tracking along the real line from $t = \rho \lambda^\alpha$ to $t=\rho\lambda^{\alpha+\xi}$; and
    \item for all $k=\alpha,\dots,\alpha+\zeta$, 
    $z(t(k))\neq z(\rho\lambda^{k})$
    where $z(t(k))$ is obtained
    by tracking along the spiral
    and $z(\rho\lambda^{k})$ is obtained
    by tracking along the real line.
\end{itemize}  
In other words, the first criterion says we have found a closed loop consisting of two pieces: a spiral arc that turns an integer number of times in the $t$-plane; and an arc that descends along the real line in~$t$. The second criterion says that if none of the points sampled along that real arc (including its endpoints) belong to the real arc descending from the starting point of the endgame, then a path interchange has occurred.
Moreover, if a path interchange occurs, then one immediately
knows that $\rho > r$. The remedy is to restart
the spiral starting from a point closer to the origin along the original real arc.
Our choice is to use $z(\rho\lambda^{\alpha+\zeta})$
obtained by tracking $z(t)$ along the real
line from $t=\rho$ to $t=\rho\lambda^{\alpha+\zeta}$
as this this point has already
been computed as part of the test for~path~closure.

\subsection{Spiral endgame}

Following the format
from the Cauchy and power
series endgames in
\S\ref{sec:Background},
the following are the steps
for spiral endgame. The method collects samples along a logarithmic spiral, $t(\tau)=\rho_0\lambda^\tau e^{2\pi i\tau}$ and also tracks one or more paths along the real ray $t\in(0,\rho_0]$. As in the Cauchy endgame, the end result is a closed path that reveals the 
purported winding number of the path along with sample points on the path that are used to estimate the endpoint. However, some extra bookkeeping and logic is necessary to handle the possibility of path interchanges.

The endgame begins at point $(z_0,t_0)=(z(\rho_0),\rho_0)$, $\rho_0\in(0,1]$. We initialize a ray $R_0=\{z_0\}$ and initialize a list of the indices of the rays crossed by the spiral as $\sJ=\{0\}$. 
\begin{description}
    \item[Sample spiral] Beginning from the last point of the spiral, track $z(t)$ around the spiral arc \mbox{$t(\tau)=\rho \lambda^\tau e^{2\pi i\tau}$} for $\tau\in[0,1]$ and record $n$ equally spaced points, $z(t(k/n))$, $k=1,\ldots,n$. 
    \item[Extend rays] For $j\in\sJ$, advance the last point of ray $R_j$ by tracking $z(t)$ along the line $t(u)=\rho\lambda^u$ for $u\in[0,1]$ and append $z(t(1))=z(\rho\lambda)$ to $R_j$.
    \item[Winding number] Check if the last point of the spiral arc matches the last point of any of the rays. If not, initialize a new ray with that endpoint, and append the index of that ray to $\sJ$. Otherwise, suppose $k$ is the index of the ray that matches. Scan backwards in $\sJ$ to find the 
    index 
    when we last crossed ray~$R_k$. If index 0 appears anywhere in the scan (including the case $k=0$), then we have a closed path consisting of, say, $\nu$ arc segments and $\nu$ segments of ray $R_k$, and that closed path includes at least one point of ray $R_0$. The closed path is declared valid and the purported winding number of the endpoint is $\nu$.
    \item[Path interchange] If the previous step detected an invalid closed path, then a path interchange has occurred. Delete the initial point of $R_0$ and restart the endgame from any of its remaining points. One option is to always restart from the final point of $R_0$. 
    \item[Estimate endpoint] If the winding number step detected a valid closed path, use the samples around the closed path to estimate the endpoint $z(0)$ using polynomial interpolation. See~\S\ref{sec:Extrapolation} for details.
    \item[Update] Set $\rho=\rho\lambda$ and loop back to the sampling step.
\end{description}

Regardless of whether the method restarts due to path interchange or updates normally, it always advances to a smaller radius in~$t$. Thus, it must eventually enter the radius of convergence of the Puiseux series wherein paths close with $\nu=\omega$. In principle, the algorithm always converges to an approximation of the path endpoint.  In practice, success depends on allocating enough computation time and enough digits of precision to keep tracking the spiral arc and the real rays to reach the desired convergence criteria.

In the \emph{Path interchange} step above, restarting from the last point of $R_0$ advances towards $t=0$ as aggressively as possible, seeking to enter the convergence radius using fewer iterations so as to eliminate costly restarts. A consequence is that if the endpoint is singular, the endgame may need to allocate more digits of precision earlier, which is also costly. There is an especially big penalty when ill-conditioning first forces an increase from double precision arithmetic performed in hardware to multiple precision arithmetic performed in software. In a situation where starting from an earlier point would avoid that penalty while the last point would incur it, it may pay to start from the earlier point, as that might be enough already to enter the convergence zone. In the experiments reported in \S\ref{sec:Examples}, we always choose the aggressive last-point option.

\subsection{Winning the spiral endgame}
As in the power series and Cauchy endgame, the spiral endgame terminates when two successive updates meet three conditions: (1)~they give the same purported winding number, (2)~their endpoint estimates agree within the specified tolerance, and (3)~the final backward error \eqref{eq:BWE} is within tolerance. The only new wrinkle is that not every update returns a purported winding number; only valid closed loops do. Accordingly, a prerequisite for meeting condition~(1) is that two successive updates give a valid closed loop.

\subsection{Pruning the spiral endgame via backtracking}\label{sec:PruningSpiral}
When the spiral endgame terminates, the index array $\sJ$ contains pointers to the rays that intersect the final closed loop. Each ray is a segment of a solution path for real $t$, one of which is the path from which the endgame was initiated. When $\omega>1$, one can backtrack along the other $\omega-1$ rays to find the starting points of these paths at $t=0$. In a homotopy for finding all solutions of a polynomial system, one can eliminate these points from the list of starting points, 
thereby eliminating the need to 
repeat the endgame for those paths.

\subsection{Branching the spiral endgame} As laid out above, when a path interchange occurs, the method discards the invalid closed loop. However, this closed loop intersects solution paths on one or more rays for real~$t$. Depending on the nature of the overall homotopy algorithm, instead of discarding these, it could be profitable to continue the spiral endgame by updating that loop to find the endpoint (or endpoints) of these alternate paths. For example, this could be part of a monodromy algorithm that uses loops in complex space to generate multiple solutions to a problem, or it could be part of a global polynomial homotopy, where these alternative rays could be backtracked to prune other starting points from the queue of work to do.

\begin{remark}
The homotopies in 
\eqref{eq:IllustrativeHomotopy}
and \eqref{eq:pathSwap}
used the gamma trick
to ensure that all solution
paths~$z(t)$ were smooth
for $t\in(0,1]$.
For example, as illustrated in~\cite[\S3.1.4]{NAGbook},
this yields circular arcs
in the complex plane
connecting the values of $1$ and $0$.  
One could replace such a circular
arc with a spiral path $t(\tau) = \lambda^\tau e^{2\pi i \tau}$ for $\tau\in[0,\infty)$
such that, for all but at most finitely many $\lambda\in(0,1)$, all solution paths
$z(t(\tau))$ are smooth.
In particular, a randomly selected $\lambda\in(0,1)$
will yield smooth solutions paths $z(t(\tau))$ 
for all $\tau\in[0,\infty)$ with probability one. Although this is a valid way to proceed, 
we prefer using the gamma trick to guarantee smoothness, leaving the selection of 
$\lambda$ free for optimizing the tradeoff between numerical conditioning and rate of convergence towards the origin.
\end{remark}

\section{Numerical Conditioning}\label{sec:condition}
In this section, we analyze the numerical conditioning of the power series, Cauchy, and spiral endgames. Each endgame is characterized by the geometric pattern of its samples in the complex $t$-plane, which through the transformation $t=s^\omega$ becomes a geometric pattern in the $s$-plane. Equations (\ref{eq:absErrorBound},\ref{eq:relErrorBound}) show that absolute error and the relative condition number for endpoint estimation are proportional to the norm of the vector of interpolation coefficients, $L(0)=\begin{bmatrix}
    \ell_{0,K}(0),\ldots,\ell_{K,K}(0)
\end{bmatrix}$. We will see that these coefficients depend on the geometric shape of the samples independent of scale. The update procedure for each endgame maintains the same shape while reducing scale. The endgames repeat this update until the samples fall within the convergence radius of the endpoint's Puiseux series and continue until the convergence criteria for winning the endgame (\S\ref{sec:Winning}) are met.

The invariance of $L$ under rescaling is a consequence of \eqref{eq:Lagrange}. Rescaling the $s$-coordinates to $s_j=\alpha\sigma_j$ gives
\[
\ell_{j,K}(0)=\prod_{\substack{j=0\\j\ne k}}^K \frac{\alpha\sigma_k}{\alpha\sigma_k-\alpha\sigma_j}=\prod_{\substack{j=0\\j\ne k}}^K\frac{\sigma_k}{\sigma_k-\sigma_j}.
\]
This shows that collecting samples in the same pattern around the origin but rescaled leaves $L(0)$ (and $||L(0)||$) unchanged. It also means that we can compute $p(0)$ by treating the data points as having been collected at $\sigma=1,\sigma_1,\ldots,\sigma_K$ with $\sigma_j=s_j/s_0$. A similar analysis shows that the Hermite coefficients, $H_{k,K}$ and $\tilde{H}_{k,K}$, are also scale invariant.

\subsection{Cauchy endgame}
The conditioning of the Cauchy endgame is clear from the estimation formula \eqref{eq:trapezoidal}. The total number of samples collected using $n$ samples per orbit and orbiting $\omega$ times to close the path is $K+1=n\omega$. Hence, each entry of $L(0)$ is $1/n\omega$ giving $||L(0)||_1=1$ and $||L(0)||_\infty=1/n\omega$. Using the 1-norm in \eqref{eq:absErrorBound}, we have 
\begin{equation}
    |\zest-z(0)| \le||\dz||_1+||R||_1.
\end{equation}
The fact that the trapezoidal rule on a uniform circular sample 
is perfectly conditioned 
makes it highly appealing for numerical work \cite{trefethen2014exponentially}.

For direct comparison to the results reported below for the power series and spiral endgames, one should use the infinity (max) norm with
\[
\max_k|\ell_{k,K}|=\max_k|H_{k,K}|=\frac{1}{K+1},\qquad 
\max_k|\tilde{H}_{k,K}|=\frac{1}{(K+1)^2}, \qquad K=n\omega-1.
\]

\subsection{Power series endgame}
The power series endgame collects samples at real values of $t$ in a geometric descent towards $t=0$. It is well-known that polynomial interpolation becomes ill-conditioned in a real interval if the points are equally spaced but stays mildly conditioned for Chebyshev spacing. This is due to the fact that Chebyshev points accumulate densely at the ends of the interval. Since the power series also accumulates points densely near the origin, we may hope that it shares the good conditioning of Chebyshev spacing. To analyze this, we must consider the rate of descent of the sampling process.

The power series samples at $t_k=\rho\lambda^k$
for $k=0,1,2,\ldots,K$
where $\lambda,\rho\in(0,1)$.
As the path approaches an endpoint with winding number $\omega$, 
the samples are transformed to $s_k=t_k^{1/\omega}=\mu q^k$
where $\mu = \rho^{1/\omega}$ and $q = \lambda^{1/\omega}$.
As $\omega$ increases
while keeping $\lambda$ and $\rho$ fixed,
both $\mu$ and $q$ approach~$1$ from below.
For what follows, it is convenient to renumber the samples in reverse, so that the closest sample to the origin is $x_0$ and the farthest is $x_K$.
Hence, we have $x_k=x_0q^{-k}$ 
for $k=0,\ldots,K$ where $0<q<1$. 

Let $p_K(x)$ be the degree $K$ polynomial that interpolates points $(x_k,z_k)$ 
for $k=0,\dots,K$.
Using the Lagrange form of polynomial interpolation
in~\eqref{eq:Lagrange}, one has
$$
    p_K(0)=\sum_{k=0}^K \ell_{k,K} z_k,
\quad\hbox{where}\quad
    \ell_{k,K}=\prod_{\substack{j=0\\j\ne k}}^K \frac{x_j}{x_j-x_k}
    = \prod_{\substack{j=0\\j\ne k}}^K \frac{1}{1-(x_k/x_j)}
    = \prod_{\substack{j=0\\j\ne k}}^K \frac{1}{1-q^{j-k}}.
$$
We can rewrite the coefficients as
\begin{align}
    \ell_{k,K}(q)
    &=\prod_{j=0}^{k-1} \frac{1}{1-q^{j-k}}\prod_{j=k+1}^{K} \frac{1}{1-q^{j-k}} \nonumber\\
    &=\prod_{j=0}^{k-1} \frac{1}{1-q^{j-k}}\prod_{j=0}^{K-k-1} \frac{1}{1-q^{j+1}}.
    \label{eq:twoProducts}
\end{align}
The two products in \eqref{eq:twoProducts} can be written in terms of the $q$-Pochhammer symbol, also known as the $q$-shifted factorial:
\begin{equation}\label{eq:qPsymbol}
(a;q)_m:=\prod_{j=0}^{m-1}(1-aq^j),\quad\text{with}\quad(a;q)_0:=1.
\end{equation}
Accordingly, we have
\begin{equation}\label{eq:ellOfqP}
    \ell_{k,K}(q) = \frac{1}{(q^{-k};q)_k}\cdot\frac{1}{(q;q)_{K-k}}=:P_k(q)Q_{K-k}(q),
\end{equation}
where the last equality defines $P_k(q)$ and $Q_{K-k}(q)$. Notice that $P_k(q)$ does not depend on~$K$.

The endgame is applied with a finite number of samples $K$. However, as $K$ grows, the trailing coefficients, that is $\ell_{k,K}(q)$ for large $k$, rapidly vanish. 
Thus, we may understand the behavior of $\ell_{k,K}(q)$ for finite $K$ by considering the limit as $K\rightarrow \infty$. The infinite product 
\begin{equation}
    (a;q)_\infty=\prod_{k=0}^{\infty}(1-aq^k)
\end{equation}
is an analytic function for $q$ inside the unit disk, i.e., $|q|<1$. The particular case of $(q;q)_\infty$ is known as Euler's function $\phi(q)$. 

\begin{figure}[!b]
    \centering
    \includegraphics[width=0.5\linewidth]{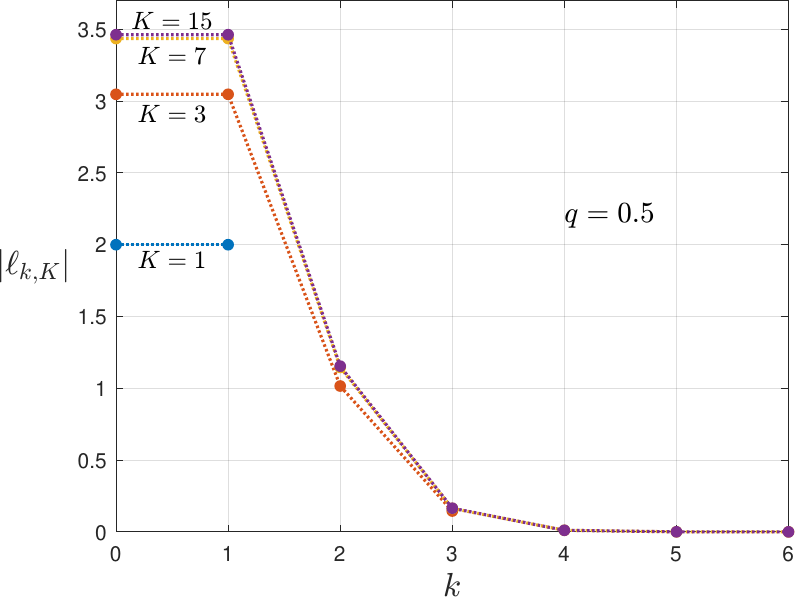}
    \caption{Absolute value of interpolation coefficients for $q=0.5$ as $K$ varies. (Values for $k>6$ are not shown as they are indistinguishable from zero at 
    this~scale.)}
    \label{fig:LshapeConvergence}
\end{figure}

For $q=0.5$, Figure~\ref{fig:LshapeConvergence} shows how the interpolation coefficients converge as the number of samples (and the order of estimation) grows. The initial coefficient $\ell_{0,K}(q)=Q_{K}(q)$ rapidly approaches~$1/\phi(q)$, which for $q=0.5$ is $\ell_{0,\infty}(q)\approx3.4627466$. 
For a fixed $k$ and $q=0.5$, this illustrates that the coefficients $|\ell_{k,K}(q)|$ 
monotonically increase as $K$ grows, but they remain finite and well-conditioned. 

We wish to understand the behavior of the infinity-norm condition number $\kappa=\max_{k=0}^n|\ell_{k,K}(q)|$ as $K$ and $q$ vary.
From \eqref{eq:ellOfqP} and the relations
\begin{align}
    (q^{-k};q)_{k}&=(1-q^{-k})(q^{-(k-1)};q)_{k-1} \label{eq:Psequence}\\
    (q;q)_{m+1}&=(1-q^{m+1})(q;q)_m, \label{eq:Qsequence}
\end{align}
one finds that
\begin{equation}\label{eq:ratio}
    \ell_{k,K}(q) = \frac{1-q^{K-k+1}}{1-q^{-k}}\ell_{k-1,K}(q):=r_k\cdot\ell_{k-1,K}(q), \quad k=1,\ldots,N.
\end{equation}
For $0<q<1$, the numerator of $r_k$ is always positive while the denominator is always negative. Thus, we see that the interpolation coefficients alternate in sign. The maximal $|\ell_{k,K}(q)|$ occurs at $k^*$ such that $|r_{k^*}|>1$ while $|r_{k^*+1}|<1$. 
For small $q$, $|r_k|<1$ for all $k$, so $|\ell_{0,K}(q)|$ is the largest coefficient. If we allow $q$ to increase, the index of the maximal coefficient increments from $k^*-1$ to $k^*$ when $|r_{k^*}|=1$. This occurs when $q=q_K^*$ 
which satisfies
$$q_K^*=\left(\frac{1+(q^*)^{K+1}}{2}\right)^{1/k^*}
\hbox{\,\,\,\,\,\,with\,\,\,\,\,\,}
q^*_\infty=\lim_{K\rightarrow\infty}q_K^*=(1/2)^{1/k^*}.
$$

Since \eqref{eq:Qsequence} implies that $Q_{K+1-k}(q)>Q_{K-k}(q)$, 
\begin{align}\label{eq:LkBound}
    \max_{0<k<K}|\ell_{k,K}(q_K^*)|<\max_{0<k<K}|\ell_{k,K+1}(q_K^*)|
    <|\ell_{k^*,\infty}(q_K^*)|=P_{k^*}(q_K^*)Q_\infty(q_K^*).
\end{align}
Table~\ref{tab:LKinfinityPower} shows the maximum interpolation coefficient for various values of $K$ and $q=q_\infty^*$ for $k^*=1,\ldots,8$. The last row is the limiting value as $K$ grows large. This shows the monotonic growth of this condition number down each column and across each row. 

\begin{table}[ht]
    \centering
\begin{tabular}{c|c|c|c|c|c|c|c|c}
    $k^*$& 1&2&3&4&5&6&7&8\\
$q_\infty^*$&0.500&0.707&0.794&0.841&0.871&0.891&0.906&0.917\\
$\max_k|\ell_{k,1}|$&2.0e+00&3.4e+00&4.8e+00&6.3e+00&7.7e+00&9.2e+00&1.1e+01&1.2e+01\\
$\max_k|\ell_{k,3}|$&3.0e+00&1.6e+01&5.0e+01&1.1e+02&2.1e+02&3.6e+02&5.7e+02&8.4e+02\\
$\max_k|\ell_{k,7}|$&3.4e+00&4.7e+01&4.2e+02&2.3e+03&9.3e+03&3.1e+04&8.9e+04&2.2e+05\\
$\max_k|\ell_{k,15}|$&3.5e+00&6.3e+01&1.3e+03&2.1e+04&2.5e+05&2.3e+06&1.6e+07&9.7e+07\\
    $|\ell_{k^*,\infty}|$&3.5e+00&6.4e+01&1.5e+03&4.1e+04&1.2e+06&3.4e+07&1.0e+09&3.2e+10\\
\end{tabular}
    \caption{Condition numbers, $\max_k|\ell_{k,K}(q)|$, as $K$ and $q$ vary for the power series endgame without derivatives (Lagrange interpolation)}
    \label{tab:LKinfinityPower}
\end{table}

In our experiments, we fix the descent factor for sampling in $t$ as $\lambda=1/2$, which means that after transformation to $s$, the descent factor is $q=(1/2)^{1/\omega}.$ With this choice, the column in Table~\ref{tab:LKinfinityPower} for $k^*=\omega$ shows the relevant condition number. One sees that for small $K$, the conditioning remains mild, but one should keep in mind that the remainder term in the Puiseux series is
\mbox{$R_K(t)=\sO(t^{(K+1)/\omega})$}. As $\omega$ grows, one would need to increase $K$ in proportion to maintain the same order of convergence.

An alternative to collecting more samples is to include the derivative at each sample using Hermite interpolation. Table~\ref{tab:HermitePower} lists the maximum absolute value of the Hermite interpolation coefficients for a range of $q$ and $K$. In comparing Tables \ref{tab:LKinfinityPower} and~\ref{tab:HermitePower}, one should keep in mind that to reach the same order of truncation error as Hermite interpolation with $K+1$ points, Lagrange interpolation must use $2K+2$ sample points. So, for example, one should compare $\max_k|H_{k,3}|$ to $\max_k|\ell_{k,7}|$. Although the Hermite conditioning is still worse when compared in this manner, it requires fewer sample points, which can be an advantage when the convergence radius is small.

\begin{table}[ht]
    \centering
\begin{tabular}{c|c|c|c|c|c|c|c|c}
$q$&0.500&0.707&0.794&0.841&0.871&0.891&0.906&0.917\\
\hline
$\max_k|H_{k,1}|$&5.0e+00&4.6e+01&1.6e+02&3.8e+02&7.4e+02&1.3e+03&2.1e+03&3.1e+03\\
$\max_k|\tilde{H}_{k,1}|$&4.0e+00&1.2e+01&2.3e+01&4.0e+01&6.0e+01&8.4e+01&1.1e+02&1.5e+02\\
\hline
$\max_k|H_{k,3}|$&1.8e+01&7.7e+02&8.5e+03&6.3e+04&3.0e+05&1.1e+06&3.1e+06&7.7e+06\\
$\max_k|\tilde{H}_{k,3}|$&9.3e+00&2.7e+02&2.5e+03&1.3e+04&4.6e+04&1.3e+05&3.2e+05&7.0e+05\\
\hline
$\max_k|H_{k,7}|$&2.6e+01&4.8e+03&5.6e+05&2.1e+07&3.7e+08&5.0e+09&5.0e+10&3.6e+11\\
$\max_k|\tilde{H}_{k,7}|$&1.2e+01&2.2e+03&1.7e+05&5.4e+06&8.6e+07&9.7e+08&7.9e+09&4.9e+10\\
\hline
$\max_k|H_{k,15}|$&2.7e+01&8.3e+03&4.3e+06&1.5e+09&2.2e+11&2.3e+13&1.3e+15&4.7e+16\\
$\max_k|\tilde{H}_{k,15}|$&1.2e+01&4.0e+03&1.6e+06&4.3e+08&6.1e+10&5.1e+12&2.6e+14&9.5e+15\\
\end{tabular}
    \caption{Condition numbers as $K$ and $q$ 
    vary for the power series endgame with derivatives. The Hermite coefficients are denoted $H_{k,K}$ and $\tilde{H}_{k,K}$ per \eqref{eq:HermiteInterpolation}.}
    \label{tab:HermitePower}
\end{table}

\subsection{Spiral endgame}\label{sec:spiralConditioning}
\begin{figure}[!b]
    \centering
    \includegraphics[width=0.5\linewidth]{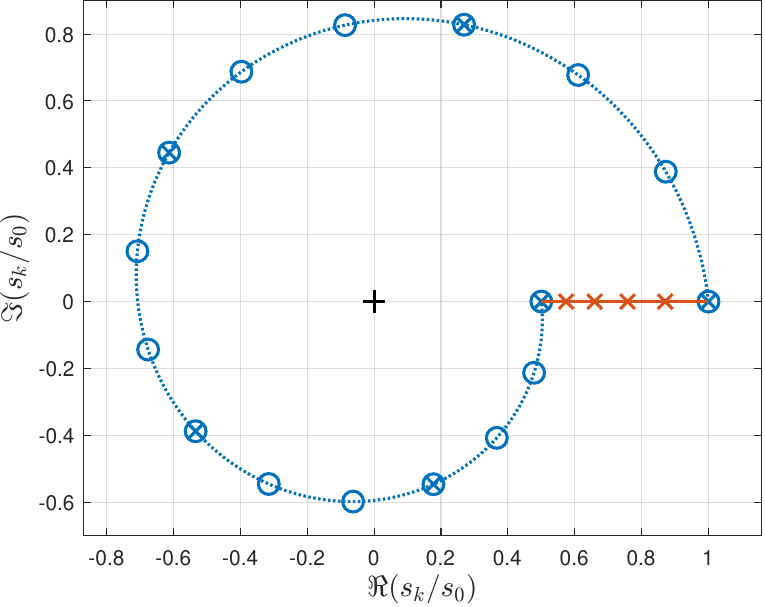}
    \caption{Sampling pattern $s_k/s_0$ for the spiral endgame with $n=3$ and $\lambda=0.5$ for an endpoint with winding number $\omega=5$. Points marked with `$\times$' correspond to having real $t$.}
    \label{fig:spiralPattern}
\end{figure}

The spiral endgame proceeds similarly to the power series endgame, except the descent factor $q$ is complex with $|q|<1$. Upon loop closure, the spiral endgame will have collected samples $s_k/s_0$, $k=0,\ldots,(n+1)\omega-1$ in the pattern shown in Figure~\ref{fig:spiralPattern}, with $n\omega+1$ points around the spiral arc and and additional $\omega-1$ points along the real line generated in testing for closure of the loop. Polynomial interpolation using just the spiral arc matches terms of the Puiseux series up to and including $t^{n\omega/\omega}=t^n$. If all the samples along the real segment are also included, this rises to $t^{n+1-1/\omega}$. If Hermite interpolation is used, these exponents double to $t^{2n}$ and $t^{2n+2-2/\omega}$, respectively.

With regard to conditioning, two types of comparison are of interest. First, we check for the impact on conditioning of including the extra real samples. They are more closely spaced than the spiral points, so conditioning will be worse.
However, if that effect is mild, the higher order of the approximation may pay off. Second, we wish to compare the spiral method's conditioning to that of the power series and Cauchy endgames.

\newcommand{\arc}{{\rm sp}}
\newcommand{\all}{{\rm all}}
For the first of these comparisons, we vary the winding number $\omega$ for $n=3$ and $\lambda=0.5$. 
Table~\ref{tab:realInVsOut} compares $\kappa_\arc=\max_k|\ell_{k,n\omega}|$ for the estimate using only points around the spiral arc 
with
\mbox{$\kappa_\all=\max_k|\ell_{k,n\omega+\omega-1}|$} for the estimate that also includes the extra samples along the real line. We see that the conditioning $\kappa_\arc$ of the spiral arc by itself remains mild for all $\omega$ in the illustrated range, but when the extra real samples are included, conditioning $\kappa_\all$ grows substantially as $\omega$ increases. This happens since
the extra real samples are all packed into the real interval $(\lambda,1)$ and, in particular, the one at $\lambda^{1-1/\omega}$ is close to~$\lambda$, the last point of the spiral arc. Due to this effect, we elect to use only the samples on the spiral arc in our experiments.

\begin{table}[!t]
    \centering
\begin{tabular}{c|c|c|c|c|c|c|c|c}
$\omega$&1&2&3&4&5&6&7&8\\
$\kappa_\arc$&9.0e-01&6.4e-01&5.5e-01&5.1e-01&5.6e-01&6.5e-01&7.4e-01&8.4e-01\\
$\kappa_\all$&9.0e-01&2.2e+00&7.2e+00&3.8e+01&2.0e+02&1.0e+03&5.0e+03&2.9e+04\\
\end{tabular}
    \caption{Conditioning of polynomial interpolation using only points on the spiral arc $(\kappa_\arc)$ versus also including points along the real line $(\kappa_\all)$, ($n=3$, $\lambda=0.5$)}
    \label{tab:realInVsOut}
\end{table}

For interpolation using only the samples on the spiral, the Lagrange coefficients take the same form as in \eqref{eq:ellOfqP} in the power series endgame with $K=n\omega$ and $q=\lambda^{1/K}e^{2\pi i/K}$. As $K$ increases, the samples are packed into the same spiral arc spaced at equal angles. The coefficients, and hence the condition number, depend only on $K$, while the order of approximation in $t$ depends only on~$n$. 

\begin{figure}[!b]
    \centering
    \includegraphics[width=0.5\linewidth]{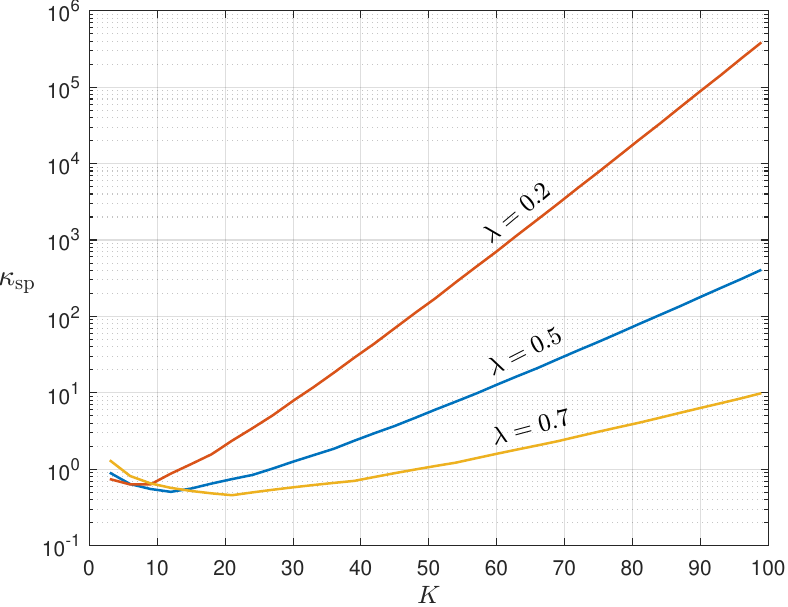}
    \caption{Condition number $\kappa_\arc$ for the spiral endgame as $K$ increases for three values of $\lambda$.}
    \label{fig:kappaArc}
\end{figure}

For $\lambda=0.2,0.5,0.7$, Figure~\ref{fig:kappaArc} shows the behavior of $\kappa_\arc$ as $K$ grows. Conditioning improves for~$\lambda$ closer to~1, where the spiral endgame approaches the Cauchy endgame. The drawback of a larger $\lambda$ is a slower descent towards the origin. A slow descent results in much wasted computation before sampling finally entering 
the endgame operating zone.
The value $\lambda=0.5$ is a compromise between the conflicting goals of descending quickly and keeping the condition number low.

As $\lambda\rightarrow1^-$, the spiral endgame becomes the Cauchy endgame, except one must be careful to drop one of the endpoints of the spiral. Otherwise, the point is repeated and the Lagrange coefficients blow up. Knocking out the repeat gives the Cauchy endgame with $\kappa_\arc=1/(K+1)=1/n\omega$.

The deterioration of conditioning shown in Figure~\ref{fig:kappaArc} as $K$ grows occurs since the sample points are distributed uniformly in angle but not in distance. The ratio formula \eqref{eq:ratio} still applies, and in particular, one sees that for $k=1$, the ratio is 
\[
  r_1=\frac{1-q^K}{1-q^{-1}}, \qquad q=\lambda^{1/K}e^{2\pi i/K}.
\]
For $\lambda=0.5$ and $K=99$, this evaluates to $|r_1|\approx 7.8$. The ratio $r_k$ decreases to approximately 1 at $k=10$, but over that interval, the growth compounds from $|\ell_{0,99}|\approx0.88$ to $|\ell_{10,99}|\approx409$. The growth could be prevented by adjusting the locations of the samples along the spiral arc, but that would destroy the self-similarity that allows the sample to be updated by advancing one orbit in $t$ to add $n$ new sample points and delete $n$ old ones. 

For a constant $\lambda$, the condition number only depends on the product $K=n\omega$. Since we use $n=3$ and $\lambda=0.5$ in all our numerical experiments, Table~\ref{tab:kappaArc} lists the corresponding condition numbers for $K=3\omega$, $\omega=1,\ldots,32$.
\begin{table}[ht]
    \centering
\begin{tabular}{c|c|c|c|c|c|c|c|c}
$K$&3&6&9&12&15&18&21&24\\
$\kappa_\arc$&9.0e-01&6.4e-01&5.5e-01&5.1e-01&5.6e-01&6.5e-01&7.4e-01&8.4e-01\\
$H$&4.5e-01&9.4e-01&1.2e+00&1.9e+00&3.0e+00&4.6e+00&6.8e+00&1.1e+01\\
$\tilde{H}$&4.1e-01&2.1e-01&1.5e-01&1.3e-01&1.6e-01&2.2e-01&2.9e-01&3.7e-01\\[2pt]
\hline 
$K$&27&30&33&36&39&42&45&48\\
$\kappa_\arc$&1.0e+00&1.3e+00&1.5e+00&1.9e+00&2.4e+00&3.0e+00&3.7e+00&4.7e+00\\
$H$&1.9e+01&3.1e+01&4.8e+01&8.5e+01&1.5e+02&2.4e+02&4.2e+02&7.3e+02\\
$\tilde{H}$&5.5e-01&8.3e-01&1.2e+00&1.8e+00&2.9e+00&4.6e+00&7.2e+00&1.2e+01\\[2pt]
\hline 
$K$&51&54&57&60&63&66&69&72\\
$\kappa_\arc$&6.1e+00&7.7e+00&9.8e+00&1.3e+01&1.6e+01&2.1e+01&2.8e+01&3.6e+01\\
$H$&1.3e+03&2.1e+03&3.8e+03&6.7e+03&1.2e+04&2.1e+04&3.7e+04&6.4e+04\\
$\tilde{H}$&1.9e+01&3.1e+01&5.1e+01&8.6e+01&1.4e+02&2.4e+02&4.0e+02&6.9e+02\\[2pt]
\hline 
$K$&75&78&81&84&87&90&93&96\\
$\kappa_\arc$&4.7e+01&6.1e+01&8.0e+01&1.0e+02&1.4e+02&1.8e+02&2.4e+02&3.1e+02\\
$H$&1.1e+05&2.0e+05&3.6e+05&6.3e+05&1.1e+06&2.0e+06&3.6e+06&6.5e+06\\
$\tilde{H}$&1.1e+03&2.0e+03&3.4e+03&5.7e+03&9.8e+03&1.7e+04&2.9e+04&5.0e+04\\[2pt]
\end{tabular}
    \caption{Conditioning of Lagrange $(\kappa_\arc)$ and Hermite $(H,\tilde{H})$  interpolation using only points on the spiral arc versus $K$ for $\lambda=0.5$. Row $H$ means $\max_k H_{k,K}$ for $H_{k,K}$ as in \eqref{eq:HermiteInterpolation}, and similarly for $\tilde{H}$.}  \label{tab:kappaArc}
\end{table}

Using Hermite interpolation to double the order of approximation gives the entries $H$ and $\tilde{H}$ in Table~\ref{tab:kappaArc}. When comparing these entries to $\kappa_\arc$, one should remember that the Hermite approximation at $K$ has the same order of approximation as Lagrange interpolation at $2K$. Taking this into account, one sees that, for the same order of approximation, the two methods have similar condition numbers. For example, Lagrange interpolation at $(n,\omega)=(6,16)$, i.e., $K=96$, has condition number $\kappa_\arc\approx3.1\cdot10^2$ which is comparable to Hermite interpolation at $(n,\omega)=(3,16)$, i.e., $K=48$, where $H\approx 7.3\cdot10^2$, and both have the same order of approximation. 
Moreover, since a first derivative is less expensive to compute than adding an extra sample point, Hermite interpolation is preferred.

\section{Experiments}\label{sec:Examples}

The Cauchy and spiral endgames are compared on several examples using the total 
number of arc segments and segments (collectively called segments) 
tracked via a parameter homotopy in Bertini~\cite{Bertini,BertiniBook} as the comparison metric. 
Due to local conditioning
and the impact of winding numbers, we do
not compare with the 
power series endgame.
Both the Cauchy and spiral endgames are initiated at $t=0.1$, use $n=3$ points along the defined contour, utilize derivatives at each sample point as described in \S\ref{Sec:Hermite}, and decrease the radius by a ratio of $\lambda=0.5$. The convergence tolerances
for winning the endgame
are $10^{-8}$
for the univariate
and bivariate examples
in \S\ref{sec:Univariate1}-\ref{sec:GO}
and $10^{-12}$
for the larger
examples in \S\ref{sec:SDP}.
In the following tables, the spiral endgame defines $t$ as the ending radius value when the 
contour has~closed.
The example files are available at \href{https://doi.org/10.7274/33940807}{doi.org/10.7274/33940807}.

\subsection{Univariate polynomial}\label{sec:Univariate1}

Returning to the univariate linear
homotopy~\eqref{eq:pathSwap}
considered in Example~\ref{ex:univariate} with $\gamma = e^{\pi i/16}$, we use the start point $(z,t) = (5,1)$ and employ both the spiral and Cauchy endgames. The spiral endgame initiated at $t=0.1$, corresponding to point 0 in Figure~\ref{fig:pathswap}, realizes a path interchange and in response re-initializes the current spiral at $t=0.0250$ along the real contour, corresponding to point 2 in Figure \ref{fig:pathswap}. With 26 segments tracked for the spiral endgame and 38 segments for Cauchy endgame, both 
endgames converge to the endpoint, to four decimal places, $-0.5000 + 0.8660i$.   The winding number and convergence metrics during the endgames are presented in Table \ref{tab:univariate_swap}.

\renewcommand{\arraystretch}{1.15}
\begin{table}[ht]
    \centering
    \begin{tabular}[2in]{c|ccc|ccc}
    \hline
    \multicolumn{1}{c}{} & \multicolumn{3}{c}{Spiral endgame} & \multicolumn{3}{c}{Cauchy endgame}\\
    \hline
    $t$ & $\nu$ & $\Vert z_j-z_{j-1}\Vert_{\infty}$ & BWE & $\nu$ & $\Vert z_j-z_{j-1}\Vert_{\infty}$ & BWE\\
    \hline
    $1.0000\cdot 10^{-1}$ & - & - & - &  $2$ & - & - \\
    $5.0000\cdot 10^{-2}$ & - & - & - &  $2$ & $1.4249\cdot 10^{-10}$ & $2.7273\cdot 10^{-1}$ \\
    $2.5000\cdot 10^{-2}$ & $\circledast$ & - & - &  $1$ & - & - \\
    $1.2500\cdot 10^{-2}$ & $1$ & - & - &  $1$ & $1.1750\cdot 10^{-2}$ & $4.8505\cdot 10^{-5}$\\
    $6.2500\cdot 10^{-3}$ & $1$ & $1.4228\cdot 10^{-3}$& $6.4785\cdot 10^{-7}$ &  $1$ & $1.1038\cdot 10^{-4}$ & $7.1298\cdot 10^{-7}$\\
    $3.1250\cdot 10^{-3}$ & $1$ & $1.4922\cdot 10^{-6}$& $1.7117\cdot 10^{-9}$ &  $1$ & $1.6210\cdot 10^{-6}$ & $1.1057\cdot 10^{-8}$\\
    $1.5625\cdot 10^{-3}$ & $1$ & $3.9397\cdot 10^{-9}$& $5.7156\cdot 10^{-12}$ &  $1$ & $2.5136\cdot 10^{-8}$ & $1.7260\cdot 10^{-10}$\\
    $7.8125\cdot 10^{-4}$ & - & - & - &  $1$ & $3.9238\cdot 10^{-10}$ & $2.6966\cdot 10^{-12}$\\
    \hline
    \end{tabular}
    \caption{Spiral and Cauchy endgames applied to the  univariate linear homotopy~\eqref{eq:pathSwap} 
    with $\gamma = e^{\pi i/16}$.
    $\circledast$ indicates the spiral endgame restarted at the corresponding $t$ value following the detection of a path interchange.
    }
    \label{tab:univariate_swap}
\end{table}

The choice of $\gamma$ can impact 
both the path 
and the occurrence of path interchanges within the spiral endgame due to how the spiral contour encircles branch points. 
Starting with $(z,t)=(5,1)$
with 
$\gamma = e^{-\pi i/16}$,
no path interchanges occur. 
Both endgames now
converge to the endpoint,
to four decimal places,
$-0.5000 - 0.8660i$.
The spiral endgame
tracked 27 segments 
compared with
38 segments for 
the Cauchy endgame. Additional details and convergence metrics available in Table~\ref{tab:univariate_noSwap}.

\begin{table}[ht]
    \centering
    \begin{tabular}[2in]{c|ccc|ccc}
    \hline
    \multicolumn{1}{c}{} & \multicolumn{3}{c}{Spiral endgame} & \multicolumn{3}{c}{Cauchy endgame}\\
    \hline
    $t$ & $\nu$ & $\Vert z_j-z_{j-1}\Vert_{\infty}$ & BWE & $\nu$ & $\Vert z_j-z_{j-1}\Vert_{\infty}$ & BWE\\
    \hline
    $1.0000\cdot 10^{-1}$ & - & - & - & $2$ & - & - \\
    $5.0000\cdot 10^{-2}$ & - & - & - & $2$ & $1.4249\cdot 10^{-10}$ & $2.7273\cdot 10^{-1}$\\
    $2.5000\cdot 10^{-2}$ & $2$ & - & - & $1$ & - & - \\
    $1.2500\cdot 10^{-2}$ & $1$ & - & - & $1$ & $1.1750\cdot 10^{-2}$ & $4.8505\cdot 10^{-5}$\\
    $6.2500\cdot 10^{-3}$ & $1$ & $1.4533\cdot 10^{-3}$ & $6.6180\cdot 10^{-7}$ & $1$ & $1.1038\cdot 10^{-4}$ & $7.1298\cdot 10^{-7}$\\
    $3.1250\cdot 10^{-3}$ & $1$ & $1.5244\cdot 10^{-6}$ & $1.7367\cdot 10^{-9}$ & $1$ & $1.6210\cdot 10^{-6}$ & $1.1057\cdot 10^{-8}$\\
    $1.5625\cdot 10^{-3}$ & $1$ & $3.9975\cdot 10^{-9}$ & $5.7637\cdot 10^{-12}$ & $1$ & $2.5136\cdot 10^{-8}$ & $1.7260\cdot 10^{-10}$\\
    $7.8125\cdot 10^{-4}$ & - & - & - & $1$ & $3.9238\cdot 10^{-10}$ & $2.6966\cdot 10^{-12}$\\
    \hline
    \end{tabular}
    \caption{Spiral and Cauchy endgames applied to the  univariate linear homotopy~\eqref{eq:pathSwap} 
    with $\gamma = e^{-\pi i/16}$
    which does not incur a path interchange.}
    \label{tab:univariate_noSwap}
\end{table}

\subsection{Univariate polynomial with multiplicity}\label{sec:Univariate2}

\begin{figure}
    \centering
    \includegraphics[width=\linewidth]{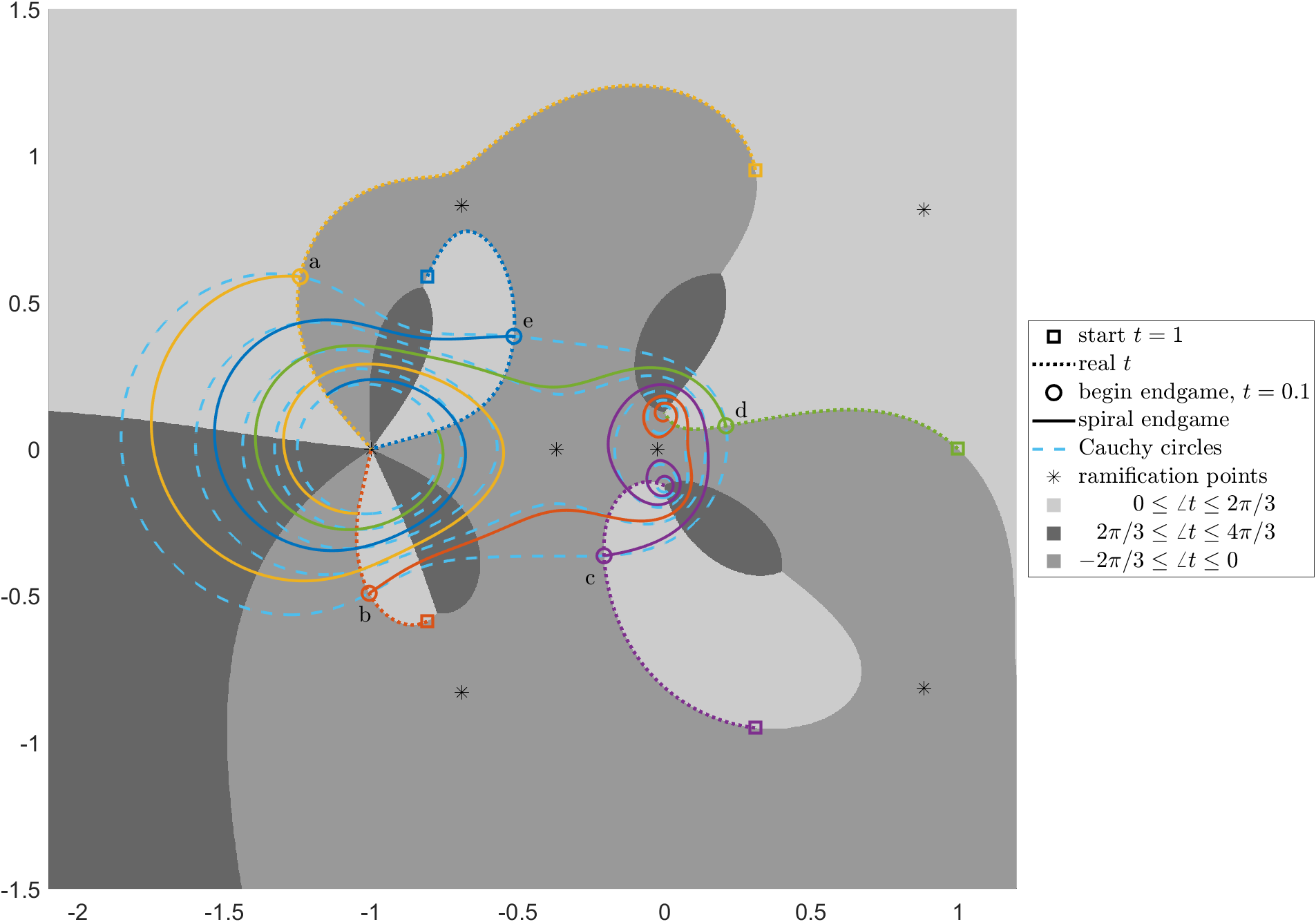}
    \caption{Univariate linear homotopy~\eqref{eq:univariateMult} with a multiplicity 3 endpoint and two nonsingular endpoints.}
    \label{fig:singularExample}
\end{figure}

Consider the linear homotopy
\begin{equation}\label{eq:univariateMult}
    h(z,t) = (1-t)(z+1)^3(z^2+1/64) + \gamma t(z^5-1)
\end{equation}
with $\gamma = e^{3i/4}$. Tracking a path from,
to four decimal places, $(z,t)=(-0.8090 - 0.5878i, 1)$, the spiral contour beginning at $t=0.1$, corresponding to point b in Figure \ref{fig:singularExample}, detects a path interchange illustrated by the red spiral converging to $z = i/8$. The current spiral is restarted along the real contour at $t=0.0125$, as outlined in \S\ref{sec:interchange}, and converges to the  endpoint $-1$ after having tracked 37 segments. In comparison, the Cauchy endgame requires a total of 82 segments.
Additional convergence metrics are available in Table~\ref{tab:univariate_mult}. 

To approximate all solutions to the target system in \eqref{eq:univariateMult}, all five paths emanating from the fifth roots of unity
are tracked. 
A total of 183 segments are utilized by the spiral endgame compared 
with the Cauchy endgame tracking 364 segments.

Pruning described in
\S\ref{sec:pruning} 
and~\S\ref{sec:PruningSpiral} can be applied to reduce computational cost.
For example, 
the data described
in Table~\ref{tab:univariate_mult} 
yields a  
purported winding number of $\nu=3$ 
which means that 
backtracking
can be used to avoid
running the same endgame
$2$ additional times.  
This is illustrated in Figure~\ref{fig:singularExample} where the spiral emanating from point 
a crosses real $t$ contours connected to points b and e before returning to the contour connected to point a. When the convergence criteria is met, a single segment can be tracked back to $t=1$ from each of the additional points along the real~$t$ contour eliminating the need 
to rerun the endgame.
This pruning reduced the total segments tracked 
by the spiral endgame to 105 and the Cauchy endgame to 202 for computing all solutions.

\begin{table}[ht]
    \centering
    \begin{tabular}[2in]{c|ccc|ccc}
    \hline
    \multicolumn{1}{c}{} & \multicolumn{3}{c}{Spiral endgame} & \multicolumn{3}{c}{Cauchy endgame}\\
    \hline
    $t$ & $\nu$ & $\Vert z_j-z_{j-1}\Vert_{\infty}$ & BWE & $\nu$ & $\Vert z_j-z_{j-1}\Vert_{\infty}$ & BWE\\
    \hline
    $1.0000\cdot 10^{-1}$ & - & - & - & $5$ & - & -\\
    $5.0000\cdot 10^{-2}$ & - & - & - & $5$ & $1.2450\cdot 10^{-7}$ & $9.4699\cdot 10^{-3}$\\
    $2.5000\cdot 10^{-2}$ & - & - & - & $3$ & - & -\\
    $1.2500\cdot 10^{-2}$ & $\circledast$ & - & - & $3$ & $1.6875\cdot 10^{-4}$ & $3.6637\cdot 10^{-18}$\\
    $6.2500\cdot 10^{-3}$ & - & - & - & $3$ & $3.1537\cdot 10^{-6}$ & $1.4894\cdot 10^{-23}$\\
    $3.1250\cdot 10^{-3}$ & - & - & - & $3$ & $5.0348\cdot 10^{-8}$ & $5.7267\cdot 10^{-29}$\\
    $1.5625\cdot 10^{-3}$ & $3$ & - & - & $3$ & $7.8879\cdot 10^{-10}$ & $2.1867\cdot 10^{-34}$\\
    $7.8125\cdot 10^{-4}$ & $3$ & $1.6872\cdot 10^{-9}$ & $2.8590\cdot 10^{-33}$ & - & - & -\\
    \hline
    \end{tabular}
    \caption{Spiral and Cauchy endgames applied to the  univariate linear homotopy~\eqref{eq:univariateMult} with $\gamma = e^{3i/4}$.    
$\circledast$ indicates the spiral endgame restarted at the corresponding $t$ value following the detection of a path interchange.
    }
    \label{tab:univariate_mult}
\end{table}   

\subsection{Griewank-Osborne}\label{sec:GO}

Consider the linear homotopy
\begin{equation}
    h(z;t) = (1-t)\left[\begin{array}{c}
29/16 z_1^3 - 2z_1z_2  \\
z_2 - z_1^2 
\end{array}\right] + \gamma t \left[\begin{array}{c}
z_1^3 - 8 \\
z_2^2 - 4
\end{array}\right]
\end{equation}
with $\gamma = 0.12 + 0.76i$
where we consider
the path 
starting from 
$(z_1, z_2, t) = (2,2,1)$.
This example highlights the increased order of approximation when including the derivatives at each sample point along the contour in the endpoint approximation, as described in \S\ref{Sec:Hermite}. When employing the spiral endgame, a path interchange is detected after the first spiral contour closes and the current spiral restarts along the real line at $t=6.2500 \cdot 10^{-3}$, indicated by $\circledast$ in Tables \ref{tab:griewankOsborne} and \ref{tab:griewankOsborne_deriv}. Without derivatives, the spiral endgame tracks 65 segments 
while the
the Cauchy endgame requires 108 segments
with the corresponding convergence metrics 
listed in Table~\ref{tab:griewankOsborne}. When derivatives are utilized in the formulation, both endgames display increased convergence rates, so they terminate earlier and track fewer segments. The spiral endgames uses 41 segments while the Cauchy endgame uses 68 segments
with additional data listed in Table~\ref{tab:griewankOsborne_deriv}. 

\begin{table}[ht]
    \centering
    \begin{tabular}[2in]{c|ccc|ccc}
    \hline
    \multicolumn{1}{c}{} & \multicolumn{3}{c}{Spiral endgame} & \multicolumn{3}{c}{Cauchy endgame}\\
    \hline
    $t$ & $\nu$ & $\Vert z_j-z_{j-1}\Vert_{\infty}$ & BWE & $\nu$ & $\Vert z_j-z_{j-1}\Vert_{\infty}$ & BWE\\
    \hline
    $1.0000\cdot 10^{-1}$ & - & - & - & $1$ & - & -\\
    $5.0000\cdot 10^{-2}$ & $1$& - & - & $1$ & $3.5592\cdot 10^{-2}$ & $2.0416\cdot 10^{-2}$\\
    $2.5000\cdot 10^{-2}$ & - & - & - & $3$ & - & -\\
    $1.2500\cdot 10^{-2}$ & - & - & - & $3$ & $6.8877\cdot 10^{-2}$ & $2.0122\cdot 10^{-3}$\\
    $6.2500\cdot 10^{-3}$ & $\circledast$ & - & - & $3$ & $7.0722\cdot 10^{-3}$ & $2.4570\cdot 10^{-4}$\\
    $3.1250\cdot 10^{-3}$ & - & - & - & $3$ & $8.6038\cdot 10^{-4}$ & $3.0621\cdot 10^{-5}$\\
    $1.5625\cdot 10^{-3}$ & - & - & - & $3$ & $1.0718\cdot 10^{-4}$ & $3.8262\cdot 10^{-6}$\\
    $7.8125\cdot 10^{-4}$ & $3$ & - & - & $3$ & $1.3392\cdot 10^{-5}$ & $4.7825\cdot 10^{-7}$\\
    $3.9063\cdot 10^{-4}$ & $3$ & $2.2518\cdot 10^{-5}$ & $5.0844\cdot 10^{-8}$ & $3$ & $1.6739\cdot 10^{-6}$ & $5.9781\cdot 10^{-8}$\\
    $1.9531\cdot 10^{-4}$ & $3$ & $3.0915\cdot 10^{-7}$ & $1.2036\cdot 10^{-8}$ & $3$ & $2.0923\cdot 10^{-7}$ & $7.4727\cdot 10^{-9}$\\
    $9.7656\cdot 10^{-5}$ & $3$ & $5.8617\cdot 10^{-8}$ & $2.6205\cdot 10^{-9}$ & $3$ & $2.6154\cdot 10^{-8}$ & $9.3408\cdot 10^{-10}$\\
    $4.8828\cdot 10^{-5}$ & $3$ & $1.0707\cdot 10^{-8}$ & $6.4556\cdot 10^{-11}$ & $3$ & $3.2693\cdot 10^{-9}$ & $1.1676\cdot 10^{-10}$\\ 
    $2.4414\cdot 10^{-5}$ & $3$ & $4.1384\cdot 10^{-10}$ & $1.0975\cdot 10^{-11}$ & - & - & -\\ 
    \hline
    \end{tabular}
    \caption{Spiral and Cauchy endgames for Griewank-Osborne without including derivatives in the endpoint approximation. $\circledast$ indicates the spiral endgame restarted at the corresponding $t$ value following the detection of a path interchange.}
    \label{tab:griewankOsborne}
\end{table}

\begin{table}[ht]
    \centering
    \begin{tabular}[2in]{c|ccc|ccc}
    \hline
    \multicolumn{1}{c}{} & \multicolumn{3}{c}{Spiral endgame} & \multicolumn{3}{c}{Cauchy endgame}\\
    \hline
    $t$ & $\nu$ & $\Vert z_j-z_{j-1}\Vert_{\infty}$ & BWE & $\nu$ & $\Vert z_j-z_{j-1}\Vert_{\infty}$ & BWE\\
    \hline
    $1.0000\cdot 10^{-1}$ & - & - & - & $1$ & - & -\\
    $5.0000\cdot 10^{-2}$ & $1$& - & - & $1$ & $7.0885\cdot 10^{-2}$ & $2.6261\cdot 10^{-2}$\\
    $2.5000\cdot 10^{-2}$ & - & - & - & $3$ & - & -\\
    $1.2500\cdot 10^{-2}$ & - & - & - & $3$ & $7.6177\cdot 10^{-2}$ & $1.1585\cdot 10^{-4}$\\
    $6.2500\cdot 10^{-3}$ & $\circledast$ & - & - & $3$ & $4.5661\cdot 10^{-4}$ & $1.6928\cdot 10^{-6}$\\
    $3.1250\cdot 10^{-3}$ & - & - & - & $3$ & $6.6661\cdot 10^{-6}$ & $2.6243\cdot 10^{-8}$\\
    $1.5625\cdot 10^{-3}$ & - & - & - & $3$ & $1.0333\cdot 10^{-7}$ & $4.0964\cdot 10^{-10}$\\
    $7.8125\cdot 10^{-4}$ & $3$ & - & - & $3$ & $1.6130\cdot 10^{-9}$ & $6.3999\cdot 10^{-12}$\\
    $3.9063\cdot 10^{-4}$ & $3$ & $5.7401\cdot 10^{-9}$ & $2.7219\cdot 10^{-12}$ & - & - & -\\
    \hline
    \end{tabular}
    \caption{Spiral and Cauchy endgames for Griewank-Osborne including derivatives in the endpoint approximation. $\circledast$ indicates the spiral endgame restarted at the corresponding $t$ value following the detection of a path interchange.}
    \label{tab:griewankOsborne_deriv}
\end{table}

\subsection{Semidefinite programming}\label{sec:SDP}

For a collection of
larger examples, we consider
computing the endpoint 
of a parameter 
homotopy arising from 
testing primal feasibility
of semi-definite programs
as follows.
We start with the 
so-called {\em messy}
instances 
for $k=2,3,4$
in the test suite of
semi-definite programs
in $\bR^{(k+1)\times(k+1)}$
that have a nonzero 
finite duality gap
created in ~\cite{pataki2018positive}.
For each of these,
we perform the 
primal feasibility 
test described in~\cite[Thm.~15]{NAGSDP}.  
This yields a single 
homotopy
path 
in~$(k+2)^2$ variables
that ends at a
feasible point
with winding number $2^{k-2}$.
Since the primal problem
is feasible but not strictly feasible,
the solution set of the homotopy
at $t=0$ is the Zariski closure
of the feasible set, resulting in an endpoint that is
not an isolated solution.
Hence, these are examples
where the endpoint is a singular
nonisolated solution
in $16$, $25$, $36$
variables for $k=2,3,4$, respectively. 
The spiral endgame exhibited one path interchange for the $k=2$ and $k=4$ cases and two path interchanges for $k=3$. The ending $t$ value and convergence metrics for the spiral endgame and Cauchy endgame are in Tables \ref{tab:SDP_spiral} and \ref{tab:SDP_cauchy}, respectively. One sees that the spiral endgame always tracks fewer segments than the Cauchy endgame, with the difference growing as the winding number increases. 
When $k=4$ with $\omega=4$, the savings is about~46\%.

\begin{table}[ht]
    \centering
    \begin{tabular}[2in]{cccccc}
    \toprule
    $k$ & $t$ & $\nu$ & $\Vert z_j-z_{j-1}\Vert_{\infty}$ & BWE & Segments\\
    \midrule
    $2$ & $1.9073\cdot 10^{-7}$ & $1$ & $1.4211\cdot 10^{-14}$ & $9.5971\cdot 10^{-17}$ & $90$\\
    $3$ & $1.9073\cdot 10^{-7}$ & $2$ & $8.8745\cdot 10^{-14}$ & $3.1670\cdot 10^{-17}$ & $94$\\
    $4$ & $1.1102\cdot 10^{-17}$ & $4$ & $6.6674\cdot 10^{-16}$ & $3.0487\cdot 10^{-17}$ & $493$\\
    \bottomrule
    \end{tabular}
    \caption{Spiral endgame convergence metrics for the semidefinite programs.}
    \label{tab:SDP_spiral}
\end{table}

\begin{table}[ht]
    \centering
    \begin{tabular}[2in]{cccccc}
    \toprule
    $k$ & $t$ & $\nu$ & $\Vert z_j-z_{j-1}\Vert_{\infty}$ & BWE & Segments\\
    \midrule
    $2$ & $9.5367\cdot 10^{-8}$ & $1$ & $1.4566\cdot 10^{-13}$ & $5.1278\cdot 10^{-17}$ & $117$\\
    $3$ & $1.9073\cdot 10^{-7}$ & $2$ & $5.6732\cdot 10^{-14}$ & $2.1274\cdot 10^{-17}$ & $134$\\
    $4$ & $7.1054\cdot 10^{-16}$ & $4$ & $3.2219\cdot 10^{-13}$ & $2.3551\cdot 10^{-17}$ & $915$\\
    \bottomrule
    \end{tabular}
    \caption{Cauchy endgame convergence metrics for the semidefinite programs.}
    \label{tab:SDP_cauchy}
\end{table}

\newpage

\section{Conclusions and Discussion}\label{sec:Conclusions}

Endgames approximate the endpoint of a homotopy path by interpolating samples along the path 
utilizing the Puiseux series expansion. This provides accurate numerical approximations even when the endpoint is singular and Newton's method may fail. Each endgame comes with an update procedure, which is iterated at exponentially decreasing radii until convergence criteria~are~met.

The proposed spiral endgame formulation offers a continuum of possibilities that includes  as its limiting cases two existing endgame methods: the power series endgame and the Cauchy endgame. The power series approach, which samples only the path of real $t$, is efficient for endpoints with low winding number, $\omega$, but suffers from poor numerical conditioning and a low order of approximation as the winding number increases. It has the most efficient update procedure, as it only requires tracking one segment to add a new point to the interpolation. In contrast, the Cauchy method has the best numerical conditioning and maintains the same order of approximation independent of the winding number, but this comes at the cost of more path segments to track. This high cost is incurred at every update, as the method initiates a new circular sample for each update, orbiting $\omega$ times around $t=0$ when in the endgame operating zone. For the same order of approximation, the condition number of the spiral endgame falls between the two extremes. Like the Cauchy method, its order of approximation is independent of winding number, but it has a less expensive update, as it orbits around $t=0$ 
just~once~per~update.

The Cauchy and spiral endgames both detect the winding number by closing a loop, a more secure approach than its counterpart in the power series endgame, which does not track loops. An additional advantage is that for $\omega>1$, the loop reveals the existence of multiple incoming paths, which in the case of numerical algebraic geometry, enables pruning of the queue of paths to be tracked. By backtracking from points where the closed loop crosses the extra incoming path, one may eliminate the cost of repeating the endgame on them.

One complication introduced by a spiral path is the phenomenon of path interchange. We detail how to detect this and correct it. The computational expense of this rectification is minimal, but it requires some extra logic in coding the method. In some cases, path interchange may even be an advantage; instead of discarding the samples of a path interchange, one may choose to continue from them to find a second solution to the system while also branching back to find the solution at the end of the original path.

Other variants of singular endgames could be considered. Once the endgame is inside the convergence radius $r$ and the correct winding number $\omega$ is known, one can use the truncated Puiseux series to predict anywhere in the disk of the $s$-plane for $|s|<r^{1/\omega}$. After using Newton's method to correct the prediction, the accurate sample (or samples) can be incorporated into a new, more accurate, approximation. This opens the question of optimizing the sample pattern for efficiency and accuracy, while also checking that the prerequisites of $r$ and $\omega$ are correct. Also, one could study whether it is always best to fit the maximal number of coefficients of the Puiseux series. An alternative would be to obtain $N$ samples and use a least-squares fit to find just $M<N$ leading coefficients instead of $N$ of them.

\section*{Acknowledgments}

\noindent JDH was supported in part by NSF CCF 2331400 and the
Robert and Sara Lumpkins Collegiate Professorship. CWW was supported in part by the Huisking
Foundation, Inc.~Research Collegiate~Professorship.

\bibliographystyle{abbrv}

\bibliography{refs,citations}

\end{document}